\documentclass{amsart}

\usepackage{a4,xy,amssymb,diagrams,times,array}
\xyoption{poly}
\xyoption{arc}
\xyoption{knot}
\xyoption{curve}
\xyoption{arrow}
\xyoption{all}

\makeindex
 
\newcommand{\N}{\mathbb{N}}
\newcommand{\C}{\mathbb{\Bbbk}}
\newcommand{\Z}{\mathbb{Z}}
\newcommand{\wis}[1]{{\text{\em \usefont{OT1}{cmtt}{m}{n} #1}}}

\newcommand{\ST}{\mathcal{S}_T}

\newcommand{\vtx}[1]{*+[o][F-]{\scriptscriptstyle #1}}

\newcommand{\centerprent}[1]{\parbox{2truecm}{\vspace{0pt} #1 }}
\newcommand{\arr}[1][\Delta]{{#1}_a}

\newcommand{\defeq}{\hspace{3.5pt}\raisebox{0.3pt}{$:$}\hspace{-3.5pt}=}
\newcommand{\Flow}{\wis{flow}}
\newcommand{\dR}{\wis{dR}}
\newcommand{\Der}{\wis{der}}
\newcommand{\map}{\longrightarrow}    
\newcommand{\set}[1]{\left\{#1\right\}}
\newcommand{\setd}[2]{\set{\ #1 \mid #2\ }}
\newcommand{\quadtext}[1]{\quad\text{#1}\quad}
\newcommand{\qquadtext}[1]{\qquad\text{#1}\qquad}     
\newcommand{\Supp}{\wis{supp}}
\newcommand{\Calo}{\wis{calo}}       
\newcommand{\imp}{\Longrightarrow}
\newcommand{\pmi}{\Longleftarrow}
\newcommand{\nin}{\not\in}
\newcommand{\Rep}{\wis{rep}}
\newcommand{\rank}{\wis{rank}}
\newcommand{\Iss}{\wis{iss}}
\newcommand{\tp}{^{\tau}}
\newcommand{\M}{M}
\newcommand{\aquot}{/\!\!/}
\newcommand{\GL}{\operatorname{GL}}

\newcommand{\PSL}{\operatorname{PSL}}
\newcommand{\gen}[1]{\langle #1 \rangle}
\newcommand{\gend}[2]{\gen{\ #1 \mid #2\ }}
\newcolumntype{C}{>{$}c<{$}}

\newtheorem{definition}{Definition}
\newtheorem{proposition}{Proposition}

\newtheorem{theorem}{Theorem}
\newtheorem{lemma}{Lemma}
\newtheorem{example}{Example}
\newtheorem{remark}{Remark}

\title{Trees of Semi-Simple Algebras}
  \author{Jan Adriaenssens}
  \address{Department of Mathematics, University of Antwerp \\ 
 Middelheimlaan 1, B-2020 Antwerp (Belgium) \\ {\tt jan.adriaenssens@ua.ac.be}}
\author{Lieven Le Bruyn} 
\address{Department of Mathematics, University of Antwerp \\ 
 Middelheimlaan 1, B-2020 Antwerp (Belgium) \\ {\tt lieven.lebruyn@ua.ac.be}}

\begin{document}
\sloppy
 
\def\ldb{\mathopen{\{\!\!\{}} \def\rdb{\mathclose{\}\!\!\}}}


 \begin{abstract}
To a tree of semi-simple algebras we associate a qurve (or formally smooth algebra) $S_T$. We introduce a Zariski- and \'etale quiver describing the finite dimensional representations of $S_T$. In particular, we show that all quotient varieties of the \'etale quiver have a natural Poisson structure induced by a double Poisson algebra structure on a certain universal localization of its path algebra.
Explicit calculations are included for the group algebras of the arithmetic groups $(P)SL_2(\Z)$ and $GL_2(\Z)$ but the methods apply as well to congruence subgroups.
\end{abstract}
\maketitle


\section{Introduction}
 
Let $B_3$ be the third braid group, generated by $\sigma_1, \sigma_2$ and satisfying the braid relation $\sigma_1 \sigma_2 \sigma_1 = \sigma_2 \sigma_1 \sigma_2$. It is well-known that $B_3$ maps surjectively onto $SL_2(\Z)$ via the map
\[
\sigma_1 \mapsto \begin{bmatrix} 1 & 1 \\ 0 & 1 \end{bmatrix} \qquad \text{and} \qquad
\sigma_2 \mapsto \begin{bmatrix} 1 & 0 \\ -1 & 1 \end{bmatrix}
\]
 Moreover, the center of $B_3$ is generated by $(\sigma_1 \sigma_2)^3$ and the quotient is isomorphic to $PSL_2(\Z) = SL_2(\Z)/ \{ \pm 1_2 \}$. Therefore, the study of finite dimensional (simple) representations of $B_3$ can be reduced to that of $SL_2(\Z)$.
 
 The group $SL_2(\Z)$ acts on the upper-half plane $\mathfrak{h} = \{ z \in \C~:~Im(z) > 0 \}$ by M\"obius transformations
\[
\begin{bmatrix} a & b \\ c & d \end{bmatrix}.z = \frac{a z + b}{c z + d} \]
 Consider the arc $y = \{ e^{i\theta}~|~\frac{\pi}{3} \leq \theta \leq \frac{\pi}{2} \}$ then it its translates under the $SL_2(\Z)$ action is an (infinite) tree
with fundamental domain the arc $\xymatrix{\vtx{Q} \ar[r]^y & \vtx{P}}$. As the stabilizer subgroups of $P$ (resp. $Q$, resp. $y$) are isomorphic to $\Z/6\Z$ (resp. $\Z/4\Z$, resp. $\Z/2\Z$) it follows from Bass-Serre theory \cite{Serre} that $SL_2(\Z)$ is the amalgamated group product
\[
SL_2(\Z) \simeq \Z/6\Z \ast_{\Z/2\Z} \Z/4\Z \]
More generally, if $G$ is a subgroup of finite index $d$ in $SL_2(\Z)$ (for example, congruence subgroups $\Gamma(N),\Gamma_0(N)$ or $\Gamma_1(N)$) then the fundamental domain $\mathfrak{h}/G$ is a Riemann surface of genus $g$ having $m$ cusps, $e_2$ elliptic points of order $2$ and $e_3$ elliptic points of order $3$ where all these numbers are related by the Riemann-Hurwitz formula
\[
g = 1 + \frac{1}{12}(d-6m-4e_2-3e_3) \]
The induced action of $G$ on the infinite tree above has a finite tree as fundamental domain and again by calculating stabilizer subgroups of vertices and edges, the groupstructure of $G$ follows as an iterated amalgamated group product, a {\em tree of groups} see \cite[\S 5]{Serre}.

As all these groups are virtually free (that is, they have a free subgroup of finite index), their groupalgebras $\C G$ are {\em qurves} (that is, quasi-free algebras in the Cuntz-Quillen terminology of \cite{CuntzQuillen} or formally smooth algebras in Kontsevich terminology \cite{KontRos}) by \cite{LBqaq}. Hence, techniques from non-commutative geometry can be used to study their finite dimensional representations. The main purpose of this paper is to initiate such study and as the main results hold in the more general setting of {\em trees of semi-simple algebras} $\ST$ they will be stated as such. We will illustrate the main results of this paper here in the special case of $SL_2(\Z)$. We should remark that in general an arithmetic group does not fit necessarily into this framework as it may contain free components. However, these components do not produce additional problems as is explained in remarks~\ref{rem1} and \ref{rem2}.

First we prove  that the study of finite dimensional representations of a tree $\ST$ of semi-simple algebras is equivalent to that of a certain universal localization of a finite dimensional hereditary algebra. For $SL_2(\Z)$ this algebra is the $22$-dimensional path algebra of the bipartite quiver $\Gamma$
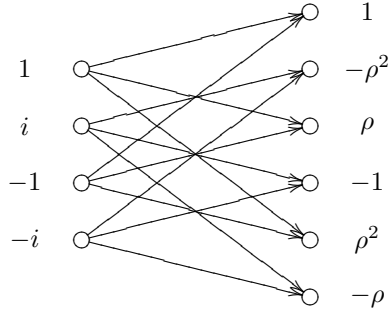
\[
\xy /r.18pc/:
\POS (40,0) *\cir<3pt>{}="b1",
 (40,10) *\cir<3pt>{}="b2",
 (40,20) *\cir<3pt>{}="b3",
 (40,30) *\cir<3pt>{}="b4",
 (40,40) *\cir<3pt>{}="b5",
 (40,50) *\cir<3pt>{}="b6",
 (0,10) *\cir<3pt>{}="a1",
 (0,20) *\cir<3pt>{}="a2",
 (0,30) *\cir<3pt>{}="a3",
 (0,40) *\cir<3pt>{}="a4",
(-10,10) *\txt{$-i$},
 (-10,20) *\txt{$-1$},
 (-10,30) *\txt{$i$},
 (-10,40) *\txt{$1$},
 (50,0) *\txt{$-\rho$},
 (50,10) *\txt{$\rho^2$},
 (50,20) *\txt{$-1$},
 (50,30) *\txt{$\rho$},
 (50,40) *\txt{$-\rho^2$},
 (50,50) *\txt{$1$}
\POS "a1" \ar "b1"
  \POS "a1" \ar "b3"
    \POS "a1" \ar "b5"
  \POS "a2" \ar "b2"
  \POS "a2" \ar "b4"
  \POS "a2" \ar "b6"
  \POS "a3" \ar "b1"
  \POS "a3" \ar "b3"
  \POS "a3" \ar "b5"
  \POS "a4" \ar "b2"
  \POS "a4" \ar "b4"
  \POS "a4" \ar "b6"
\endxy
\]
($\rho$ is a primitive $3$rd root of unity). The fact that $\Gamma$ is bipartite corresponds to the fact that the fundamental domain-tree for $SL_2(\Z)$ is reduced to one edge $\xymatrix{\vtx{Q} \ar[r]^y & \vtx{P}}$.
The left vertices of $\Gamma$ correspond to the simple representation of the stabilizer subgroup of $Q$ (that is, $\Z/4\Z$) and the right vertices to those of $\Z/6\Z$ (being the stabilizer subgroup of $P$). An arrow is drawn from a left to a right vertex whenever there is a $\Z/2\Z$-morphism between the simples ($\Z/2\Z$ being the stabilizer subgroup of the edge). 

To an $n$-dimensional representation $V$ of $SL_2(\Z) \simeq \Z/4\Z \ast_{\Z/2\Z} \Z/6\Z$ we associate a representation of the quiver $\Gamma$ where the left (resp. right) vertex-spaces give a decomposition of $V \downarrow_{\Z/4\Z}$ (resp. of $V \downarrow_{\Z/6\Z}$) into simples and where the arrows are the components of the $\Z/2\Z$-isomorphism 
\[
(V \downarrow_{\Z/4\Z}) \downarrow_{\Z/2\Z} \rTo^{\sim} (V \downarrow_{\Z/6\Z}) \downarrow_{\Z/2\Z}
\]
That is, $V$ determines an $n$-dimensional representation of the universal localization $\C \Gamma_{\Sigma}$ where $\Sigma$ is the morphism between projective right $\C \Gamma$-modules
determined by the matrix
\[
\Sigma = \begin{bmatrix}
a_{11} & 0 & a_{31} & 0 \\
0 & a_{22} & 0 & a_{42} \\
a_{13} & 0 & a_{33} & 0 \\
0 & a_{24} & 0 & a_{44} \\
a_{15} & 0 & a_{35} & 0 \\
0 & a_{26} & 0 & a_{46}
\end{bmatrix}
\]
where $a_{ij}$ is the arrow from the $i$-th left vertex to the $j$-th right vertex. Alternatively, finite dimensional $SL_2(\Z)$-representations determine $\theta$-semistable $\Gamma$-representations (where $\theta = (-1,-1,-1,-1;1,1,1,1,1,1)$) for which the matrix $\Sigma$ becomes invertible. As such the quotient varieties $\wis{iss}_n~SL_2(\Z)$ representing isomorphism classes of semi-simple $n$-dimensional $SL_2(\Z)$-representations decompose into irreducible components $\wis{iss}_{\alpha}~SL_2(\Z)$ which are affine open pieces of moduli spaces  $\wis{moduli}^{\theta}_{\alpha}~\Gamma$ of $\theta$-semistable representations of $\Gamma$ as introduced in \cite{King}.

It turns out that moduli spaces of quiver representations are much harder to study than algebraic quotient varieties of quiver representations. In our previous paper \cite{JALB1} we reduced local quivers to reduce the \'etale local structure of these moduli spaces to quiver quotient varieties. From \cite{LBqaq} we recall that all these local quiver settings are determined by one {\em \'etale quiver} $\Psi$ associated to $SL_2(\Z)$. The vertices of $\Psi$ correspond to those components $\wis{iss}_{\alpha}~SL_2(\Z)$ consisting entirely of simple representations. These are precisely the $12$ one-dimensional representations $\{ S_1,\hdots,S_{12} \}$ of $SL_2(\Z)$ corresponding to the following dimension vectors for $\Gamma$
\[
\begin{cases}
 g_1 &=(1,0,0,0,1,0,0,0,0,0) \\
 g_2 &=(0,0,1,0,0,0,1,0,0,0) \\
 g_3 &=(1,0,0,0,0,0,0,0,1,0) \\
 g_4 &=(0,0,1,0,1,0,0,0,0,0) \\
 g_5 &=(1,0,0,0,0,0,1,0,0,0) \\
 g_6 &=(0,0,1,0,0,0,0,0,1,0) 
\end{cases} \qquad
\begin{cases}
 g_7 &=(0,1,0,0,0,1,0,0,0,0) \\
 g_8 &=(0,0,0,1,0,0,0,1,0,0) \\
 g_9 &=(0,1,0,0,0,0,0,1,0,0) \\
 g_{10} &= (0,0,0,1,0,1,0,0,0,0) \\
 g_{11} &= (0,1,0,0,0,0,0,1,0,0) \\
 g_{12} &=(0,0,0,1,0,0,0,0,0,1)
\end{cases}
\]
The number of arrows between two vertices in $\Psi$ is given by $-\chi_{\Gamma}(g_i,g_j)$ where $\chi_{\Gamma}$ is the Euler-form of $\Gamma$. The \'etale quiver $\Psi$ for $SL_2(\Z)$ has the following form
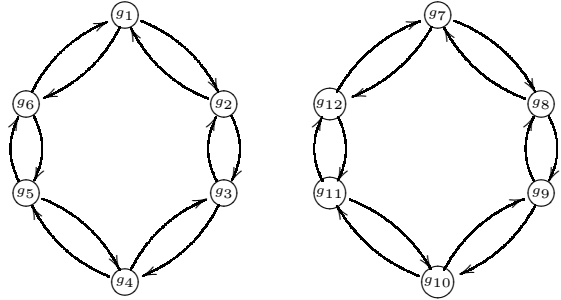
\[
\xymatrix{& \vtx{g_1} \ar@/^/[ld] \ar@/^/[rd] & \\
\vtx{g_6} \ar@/^/[ru] \ar@/^/[d] & & \vtx{g_2} \ar@/^/[lu] \ar@/^/[d] \\
\vtx{g_5} \ar@/^/[u] \ar@/^/[rd] & & \vtx{g_3} \ar@/^/[u] \ar@/^/[ld] \\
& \vtx{g_4} \ar@/^/[lu] \ar@/^/[ru] &}~\qquad~\qquad
\xymatrix{
& \vtx{g_7} \ar@/^/[ld] \ar@/^/[rd] & \\
\vtx{g_{12}} \ar@/^/[ru] \ar@/^/[d] & & \vtx{g_8} \ar@/^/[lu] \ar@/^/[d] \\
\vtx{g_{11}} \ar@/^/[u] \ar@/^/[rd] & & \vtx{g_9} \ar@/^/[u] \ar@/^/[ld] \\
& \vtx{g_{10}} \ar@/^/[lu] \ar@/^/[ru] & }
\]
The quiver $\Psi$ can be used to determine all components $\wis{iss}_{\alpha}~SL_2(\Z)$ which contain a Zariski open subset of simple representations. Write
\[
\alpha = a_1 g_1 + \hdots + a_6 g_6 + b_1 g_7 + \hdots + b_6 g_{12} \]
then this is the case (using the results from \cite{LBProcesi}) if and only if the following inequalities are satisfied for all $1 \leq i \leq 6$
\[
a_i \leq a_{i+1} + a_{i-1} \quad \text{and} \quad b_i \leq b_{i+1}+b_{i-1} \]
where indices are taken modulo $6$. Another application of $\Psi$ is to determine the dimension vectors $\alpha$ such that the quotient variety $\wis{iss}_{\alpha}~SL_2(\Z)$ is smooth. Using the results from \cite{Bocklandt1} this is the case if and only if $\alpha$ corresponds to the dimension vectors $\vec{a} = (a_1,\hdots,a_6)$ and $\vec{b} = (b_1,\hdots,b_6)$ of $\Psi$ where both components are one of the following (upto cyclic permutation of vertices)
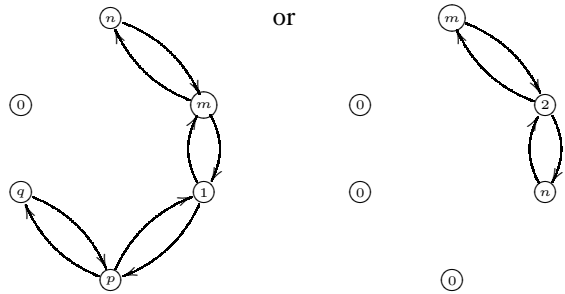
\[
\xymatrix{ & \vtx{n}  \ar@/^/[rd] & \\ \vtx{0} & & \vtx{m} \ar@/^/[lu] \ar@/^/[d] \\ \vtx{q}  \ar@/^/[rd] & & \vtx{1} \ar@/^/[u] \ar@/^/[ld] \\ & \vtx{p} \ar@/^/[lu] \ar@/^/[ru] & } \qquad \text{or} \qquad
\xymatrix{& \vtx{m} \ar@/^/[rd] & \\
\vtx{0} & & \vtx{2} \ar@/^/[lu] \ar@/^/[d] \\
\vtx{0} & & \vtx{n} \ar@/^/[u] \\
& \vtx{0} & }
\]
for some $m,n,p,q \in \N$. This is a consequence of the general fact that there is an \'etale isomorphism between the quotient variety $\wis{iss}_{(\vec{a},\vec{b})}~\Psi$ near the point corresponding to the zero representation and the quotient variety $\wis{iss}_{\alpha}~SL_2(\Z)$ near the semi-simple representation
\[
M_{\alpha} = S_1^{\oplus a_1} \oplus \hdots \oplus S_6^{\oplus a_6} \oplus S_7^{\oplus b_1} \oplus \hdots \oplus S_{12}^{\oplus b_6} \]
Hence, studying the quiver quotient variety $\wis{iss}_{(\vec{a},\vec{b})}~\Psi$ near $\overline{0}$ allows us to construct a Zariski open subset of $\wis{iss}_{\alpha}~SL_2(\Z)$ near $M_{\alpha}$.
More generally, the \'etale quiver $\Psi$ contains enough information to determine the local quiver setting whose quotient variety gives the \'etale local structure of $\wis{iss}_{\alpha}~SL_2(\Z)$ in {\em any} point, see \cite{LBqaq} for more details.

Observe that the \'etale quiver $\Psi$ of $SL_2(\Z)$ is the double quiver of two extended Dynkin quivers of type $\tilde{A}_5$ and as such carries a non-commutative symplectic structure, see \cite{Ginzburg} or \cite{LBBocklandt}. As a consequence, all quotient varieties $\wis{iss}_{\alpha}~SL_2(\Z)$ are locally symplectic varieties. For an arbitrary tree $\ST$ of semi-simple algebras this is no longer the case. However, we will prove that in general the quotient varieties $\wis{iss}_{\alpha}~\ST$ are locally Poisson varieties, where the Poisson structure is induced from a double Poisson structure (introduced in \cite{VdBPoisson}) on a certain universal localization of the path algebra of a double quiver determined by the \'etale quiver of $\ST$. 

The \'etale quiver $\Psi$ corresponding to such a tree $\ST$ of semi-simple algebras is shown to be always symmetric, and when there is an even amount of loops in each vertex it will be even double. This allows us to define the necklace Lie algebra and with it vectorfields and flows on the quiver quotient varieties. Indeed, one-way necklaces will correspond to locally nilpotent symplectic derivations, giving rise to symplectic flows on the path algebra and on the affine quotient varieties $\wis{iss}_\alpha~\ST$. For $SL_2(\Z)$, this allows for a very straightforward description but for $GL_2(\Z)$ the construction is slightly more exotic. Finally we turn these flows in to integral curves by constructing natural compactifications of these quotient varieties. 

The present paper is an updated and expanded version of the unpublished paper \cite{AdriLBvorig}.

\section{The Zariski quiver} \label{Zariski}
 
 In this section we introduce the main topic of investigation : {\em trees of semi-simple algebras}. We will show that the representation varieties of these {\em qurves} can be identified with those of certain universal localization of finite dimensional hereditary algebras. As such, all these varieties are rational varieties.
 
 Throughout, we work over an algebraically closed field $\C$ of characteristic zero. Recall that if $S$ is a (finite dimensional) semi-simple $\C$-algebra and if $A,B$ are $S$-algebras, then the {\em amalgamated free product} $A \ast_S B$ is the $\C$-algebra determined by the following universal property. For all $\C$-algebras $D$ and all $\C$-algebra morphisms $\phi~:~A \rTo D$ and $\psi~:~B \rTo D$ such that $\phi | S = \psi | S$ there is a uniquely determined $\C$-algebra morphism $A \ast_S B \rTo D$ making the diagram below commutative.
\[
\xymatrix{&  A \ar[rd]^{\phi} \ar[d] & \\
 S \ar[ru] \ar[rd] & A \ast_S B \ar@{.>}[r] & D \\
 & B \ar[ru]^{\psi} \ar[u] &}
\]
 Recall from \cite{LBqaq} that a {\em qurve} is a $\C$-algebra $A$ such that the universal bimodule $\Omega^1_{\C}(A)$ of derivations is a projective $A$-bimodule. Remark that {\em qurves} are called {\em quasi-free algebras} in \cite{CuntzQuillen} and {\em formally smooth algebras} in \cite{KontRos}. If $A$ and $B$ are qurves, then so is the amalgamated product $A \ast_S B$ by \cite[lemma 1]{LBqaq}.
 
\begin{definition} A {\em tree of semi-simple algebras} $\mathcal{S}_T$ is determined by the following data
\begin{itemize}
\item{A finite tree $T$ with vertices $V$ and edges $E$.}
\item{For every vertex $v \in V$ a semi-simple $\C$-algebra $S_v$.}
\item{For every edge $\xymatrix{\vtx{v} \ar@{-}[r]^e & \vtx{w}} \in E$ a semi-simple $\C$-algebra $S_e$ and distinguished embeddings $S_e \rInto^{i_{e,v}} S_v$ and $S_e \rInto^{i_{e,w}} S_w$.}
\end{itemize}
 The algebra $\mathcal{S}_T$ is the $\C$-algebra determined by the universal property that for every $\C$-algebra $B$ and $\C$-algebra morphisms $\phi_v~:~S_v \rTo B$ satisfying $\phi_v | S_e = \phi_w | S_e$ for every edge $\xymatrix{\vtx{v} \ar@{-}[r]^e & \vtx{w}} \in E$ there is a unique $\C$-algebra morphism $\mathcal{S}_T \rTo B$ making the above diagrams for every edge commutative.
\end{definition}
 
 The algebra $\mathcal{S}_T$ can be constructed inductively by considering dropping an edge $\xymatrix{\vtx{v} \ar@{-}[r]^e & \vtx{w}}$ such that $w$ is a leaf of $T$ and replacing $S_v$ by the amalgamated product $S_v \ast_{S_e} S_w$ and iterating this process until the tree is reduced to its root. As all algebras obtained satisfy the required universal property, $\mathcal{S}_T$ does not depend on the choices made. From this construction the next result follows.
 
\begin{proposition} A tree of semi-simple algebras $\mathcal{S}_T$ is a qurve.
\end{proposition}
 
 From this and \cite{LBnagatn} it follows that for any natural number $n \in \N$ the affine scheme $\wis{rep}_n~\mathcal{S}_T$ of $n$-dimensional representations of $\mathcal{S}_T$ is a smooth variety, though it may have several connected (= irreducible) components
\[
\wis{rep}_n~\ST = \bigsqcup_{| \alpha | = n} \wis{rep}_{\alpha}~\ST \]
 The label $\alpha$ determining such a component is called a {\em dimension vector} of total dimension $n$. We will show that these dimension vectors $\alpha$ actually correspond to certain dimension vectors of a finite quiver $\Gamma = \Gamma(\ST)$.
 
 For a tree $\ST$ of semi-simple algebras we have decompositions of the vertex- and edge-algebras in matrix-algebras
\[
 S_v = M_{d_v(1)}(\Bbbk) \oplus \hdots \oplus M_{d_v(n_v)}(\Bbbk) \quad \text{resp.} \quad
 S_e = M_{d_e(1)}(\Bbbk) \oplus \hdots \oplus M_{d_e(n_e)}(\Bbbk) \]
 The embedding $S_e \rInto S_v$ can be depicted by a Bratelli-diagram or, equivalently, by natural numbers $a^{(ev)}_{ij}$ for $1 \leq i \leq n_e$ and $1 \leq j \leq n_v$ satisfying the numerical restrictions
\[
 d_v(j) = \sum_{i=1}^{n_e} a^{(ev)}_{ij} d_e(i) \qquad \text{for all $1 \leq j \leq n_v$ and all $v \in V$ and $e \in E$} \]
 Remark that these numbers give the {\em restriction data}, that is, the multiplicities of the simple components of $S_e$ occurring in the restriction $V^{(v)}_j \downarrow_{S_e}$ for the simple components $V_j^{(v)}$ of $S_v$. 
 
 From these decompositions and Schur's lemma it follows that for any edge
 $\xymatrix{\vtx{v} \ar@{-}[r]^e & \vtx{w} }$ in the tree $T$ we have
\[
 Hom_{S_e}(V_i^{(v)},V_j^{(w)}) = \sum_{k=1}^{n_e} a_{ki}^{(ev)} a_{kj}^{(ew)} = n_{ij}^{(e)} \]
 
\begin{definition}
 For a tree $\ST$ of semi-simple algebras we define the quiver
 $\Gamma = \Gamma(\ST)$ as follows 
\begin{itemize}
\item{Vertices :  for any vertex $v \in V$ of $T$ take $n_v$ vertices $\{ \mu^{(v)}_1,\hdots,\mu^{(v)}_{n_v} \}$ (corresponding to the matrix-components of $S_v$). }
\item{Arrows : take the unique root-to-leaves orientation $\vec{T}$ on all the edges of $T$. For an edge $\xymatrix{\vtx{v} \ar@{-}[r]^e & \vtx{w} }$ in $T$ we add for each $1 \leq i \leq n_v$ and each $1 \leq j \leq n_w$ precisely $n^{(e)}_{ij}$ arrows between the vertices $\mu^{(v)}_i$ and $\mu^{(w)}_j$ oriented in the same way as the edge $e$ in $\vec{T}$.}
\end{itemize}
 We call $\Gamma$ the {\em Zariski quiver} of the tree of semi-simple algebras $\ST$.
\end{definition}

\begin{example} \label{exgl2}
The integral general linear group $GL_2(\Z)$ is the amalgamated free group product of dihedral groups $D_6 \ast_{D_2} D_4$ (see for example \cite[Thm. 23.1]{Zie}). Consider a hexagon (having symmetry group $D_6$) and a square (having symmetry group $D_4$) lined up as
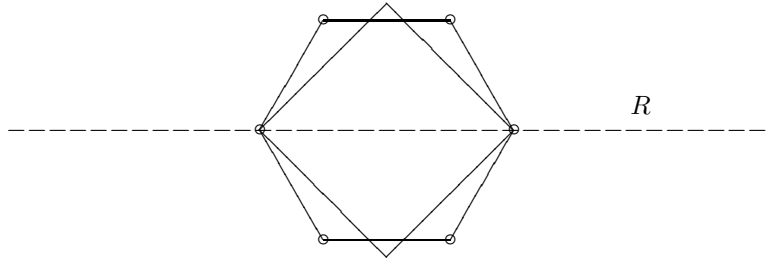
\[
\xy/r4pc/:{\xypolygon6{\circ}},
{\xypolygon4{~={90}}},
(-2,0);
(4,0) **{-};
(3,.2) *{R} \endxy
\]
$D_4$ is generated by $U$ (rotation by $90^o$) and $R$ (symmetry along the line) and has $5$ conjugacy classes : $\{ 1,R,R',U,C \}$ where $C$ is the central point reflection and $R'$ is reflection along a line through two midpoints of edges of the square. $D_6$ is generated by $V$ (rotation over $60^o$) and $R$ and has $6$ conjugacy classes " $\{ 1,R,R",V,V^2,C \}$ where $R"$ is a reflection along a line through two midpoints of the hexagon. The common subgroup $D_2$ is generated by $C$ and $R$ and is Klein's Vierergruppe. Let $\{ X_1,\hdots,X_5 \}$ (resp. $\{ Y_1,\hdots,Y_6 \}$ and $\{ Z_1,\hdots,Z_4 \}$ be the simple representations of $D_4$ (resp. of $D_6$ and $D_2$). We obtain the restriction data  
\[
\begin{array}{c|cccc}
D_4 \downarrow_{D_2} & 1 & C & R & R \\
\hline
X_1 \downarrow_{D_2} & 1 & 1 & 1 & 1 \\
X_2  \downarrow_{D_2}  & 1 & 1 & -1 & -1 \\
X_3  \downarrow_{D_2}  & 1 & 1 & 1 & 1 \\
X_4  \downarrow_{D_2}  & 1 & 1 & -1 & -1 \\
X_5  \downarrow_{D_2}  & 2 & -2 & 0 & 0
\end{array}~\quad \text{and}~\quad \begin{array}{c|cccc}
D_6 \downarrow_{D_2} & 1 & C & R & R" \\
\hline
Y_1 \downarrow_{D_2}  & 1 & 1 & 1 & 1 \\
Y_2  \downarrow_{D_2} & 1 & 1 & -1 & -1 \\
Y_3 \downarrow_{D_2} & 1 & -1 & -1 & 1 \\
Y_4 \downarrow_{D_2} & 1 & -1 & 1 & -1 \\
Y_5 \downarrow_{D_2} & 2 & -2 & 0 & 0 \\
Y_6 \downarrow_{D_2} & 2 & 2 & 0 & 0
\end{array}
\]
whence
\[
X_1 \downarrow_{D_2} = Z_1 \quad X_2 \downarrow_{D_2} = Z_4 \quad X_3 \downarrow_{D_2} = Z_1 \quad X_4 \downarrow_{D_2} =Z_4 \quad X_5 \downarrow_{D_2} =Z_2 \oplus Z_3 \]
\[
Y_1 \downarrow_{D_2} =Z_1 \quad Y_2 \downarrow_{D_2} =Z_4 \quad Y_3 \downarrow_{D_2} =Z_2 \]
\[
Y_4 \downarrow_{D_2} =Z_3 \quad Y_5 \downarrow_{D_2} =Z_2 \oplus Z_3 \quad Y_6 \downarrow_{D_2} = Z_1 \oplus Z_4 \]
Therefore, the Zariski quiver $\Gamma(GL_2(\Z))$ has the following form
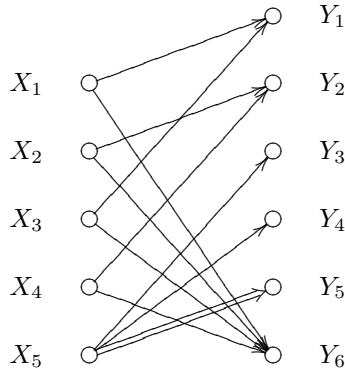
\[
\xymatrix@=.4cm{
& & & & & \vtx{} & Y_1 \\
X_1 & \vtx{} \ar[rrrru] \ar[rrrrdddd] & & & & \vtx{} & Y_2 \\
X_2 & \vtx{} \ar[rrrru] \ar[rrrrddd] & & & & \vtx{} & Y_3 \\
X_3 & \vtx{} \ar[rrrruuu] \ar[rrrrdd] & & & & \vtx{} & Y_4 \\
X_4 & \vtx{} \ar[rrrruuu] \ar[rrrrd] & & & & \vtx{} & Y_5 \\
X_5 & \vtx{} \ar[rrrruuu] \ar[rrrruu] \ar@{=>}[rrrru] & & & & \vtx{}  & Y_6 }
\]
which decomposes into the two connected components given in the introduction.
\end{example}

Observe that since $T$ is a tree, there are no oriented cycles in the Zariski quiver $\Gamma = \Gamma(\ST)$ of a tree $\ST$ of semi-simple algebras. As a consequence, the {\em path algebra} $\C \Gamma$ is a finite dimensional hereditary algebra. One can show, see \cite[Prop. 8]{LBqaq}, that any $n$-dimensional representation of $\ST$ is isomorphic to an $\alpha$-dimensional representation of $\C \Gamma$ where the dimension vector $\alpha=(\alpha_i^{(v)}~:~v \in V, 1 \leq i \leq n_v)$ is an {\em $n$-dimension vector}, that is satisfies the following numerical restrictions 
\begin{itemize}
\item{For every vertex $v \in V$ we have $
\sum_{i=1}^{n_v} d_v(i) \alpha^{(v)}_i = n$.}
\item{For every edge  $\xymatrix{\vtx{v} \ar@{-}[r]^e & \vtx{w} }$ and all $1 \leq j \leq n_e$ we have
\[
\sum_{i=1}^{n_v} \alpha_i^{(v)} a_{ji}^{(ev)} = \sum_{k=1}^{n_w} \alpha_k^{(w)} a_{jk}^{(ew)} \]}
\end{itemize}
These restrictions encode the requirement that any $n$-dimensional representation $V$ of $\ST$ when restricted to a vertex semi-simple algebra $S_v$ decomposes into a direct sum of simples and that the identification of two such decompositions $V \downarrow S_v$ and $V \downarrow S_w$ for an edge
$\xymatrix{\vtx{v} \ar@{-}[r]^e & \vtx{w} }$ must be an $S_e$-homomorphism. We refer to \cite[\S 7]{LBqaq} for more details. An alternative proof uses a certain {\em universal localization}  $\C \Gamma_{\Sigma}$ of the finite dimensional algebra $\C \Gamma$. 

Recall that every vertex $\mu_i^{(v)}$ of $\Gamma(\ST)$ determines an indecomposable projective right $\C \Gamma$ module
\[
P_{v,i} = e_{\mu^{(v)}_i} \C \Gamma \]
where $e_{\mu^{(v)}_i}$ is the vertex-idempotent. Observe that a $\C$-basis of $P_{v,i}$ is given by the paths in $\Gamma$ terminating in the vertex $\mu_i^{(v)}$. Hence, if $\xymatrix{\vtx{v} \ar@{-}[r]^e & \vtx{w} }$ is an edge in the tree $T$ and if $\xymatrix{\vtx{\mu_i^{(v)}} \ar@{->}[r]^a & \vtx{\mu_j^{(w)}} }$ is an arrow in $\Gamma(\ST)$, then left-multiplication by $a$ determines a right $\C \gamma$-module morphism
$P_{i,v} \rTo^{a.} P_{j,w}$ between the indecomposable projective right $\C \Gamma$-modules. In fact, for the given edge $\xymatrix{\vtx{v} \ar@{-}[r]^e & \vtx{w} }$ we will define a right $\C \Gamma$-module morphism
\[
P_{1,v}^{\oplus \sum_{i=1}^{n_e} a_{i1}^{(ev)}} \oplus \hdots \oplus P_{n_v,v}^{\oplus \sum_{i=1}^{n_e} a_{in_v}^{(ev)} } \rTo^{\sigma_e} P_{1,w}^{\oplus \sum_{i=1}^{n_e} a_{i1}^{(ew)}} \oplus \hdots \oplus P_{n_w,w}^{\oplus \sum_{i=1}^{n_e} a_{in_w}^{(ew)}} \]
determined by the $\sum_{j=1}^{n_w} \sum_{i=1}^{n_e} a_{ij}^{(ew)} \times \sum_{j=1}^{n_v} \sum_{i=1}^{n_e} a_{ij}^{(ev)}$ matrix with entries zeroes except for exactly $\sum_{i=1}^{n_v} \sum_{j=1}^{n_w} n_{ij}^{(e)}$ entries (each corresponding to one of the paths in $\Gamma$ corresponding to the edge $e$) at the appropriate places to encode the fact that these paths describe an $S_e$-morphism.

\begin{example} One of the two components of $\Gamma(GL_2(\Z))$ corresponds to the situation
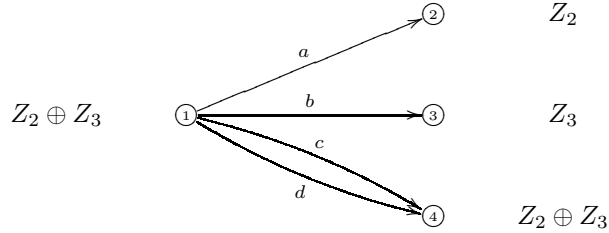
\[
\xymatrix{
& & & & \vtx{2} & Z_2 \\
Z_2 \oplus Z_3 & \vtx{1} \ar[rrru]^a \ar[rrr]^b \ar@/^1ex/[rrrd]^c  \ar@/_1ex/[rrrd]_d& & & \vtx{3} & Z_3 \\
& & & & \vtx{4} & Z_2 \oplus Z_3 } \]
whence
\[
\begin{bmatrix}
a & 0 \\
0 & b \\
c & 0 \\
0 & d \end{bmatrix}~:~P_{1}^{\oplus 2} \rTo P_{2} \oplus P_{3} \oplus P_{4}^{\oplus 2} \]
is the corresponding morphism of right $\C \Gamma$-modules.
\end{example}

Repeating this procedure for every edge $e \in E$ we obtain a set $\Sigma = \{ \sigma_e~:~e \in E \}$ of right $\C \Gamma$-module morphisms. Recall from \cite{Schofield} that the {\em universal localization} $\C \Gamma_{\Sigma}$ is the algebra with the property that all extended module morphisms $\sigma_e \otimes \C \Gamma_{\Sigma}$ are invertible and 'minimal' as such, that is, whenever there is a $\C$-algebra morphism $\C \Gamma \rTo B$ such that all $\sigma \otimes B$ are invertible, then there is an algebra morphism $\C \Gamma_{\Sigma} \rTo B$ inducing these isomorphisms.

\begin{theorem} For a tree $\ST$ of semi-simple algebras, there is a natural one-to-one correspondence between
\begin{itemize}
\item{Isomorphism classes of finite dimensional $\ST$-representations.}
\item{Isomorphism classes of finite dimensional $\C \Gamma(\ST)_{\Sigma}$-representations.}
\end{itemize}
\end{theorem}

\begin{proof} Left to the reader. Compare with \cite[Prop 8]{LBqaq}.
\end{proof}

\begin{remark} \label{rem1}
We can extend the setting to allow for free group or free algebra components. If the vertex algebra $S_v \simeq \C[x]$ or $\C[x,x^{-1}]$ is connected to the tree by $\xymatrix{\vtx{v} \ar[r]^e & \vtx{w}}$ then $S_e \simeq \C$ and the corresponding bipartite part of the Zariski quiver has the form
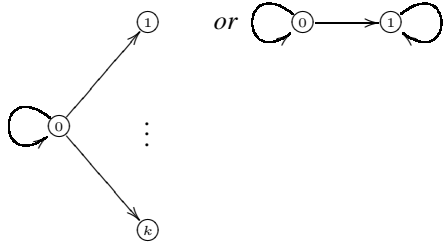
\[
\xymatrix{&  \vtx{1} \\
\vtx{0} \ar@(ul,dl) \ar[ru] \ar[rd] & \vdots \\
& \vtx{k}} \qquad \text{or} \qquad \xymatrix{\vtx{0} \ar@(ul,dl) \ar[r] & \vtx{1} \ar@(ur,dr)} \]
depending on whether $S_w$ is a semi-simple algebra on $k$ components or is $\simeq \C[x]$ or $\C[x,x^{-1}]$. The rest of the arguments can then be repeated verbatim.
\end{remark}

\section{The Etale quiver}
 
 The results of the previous section can be summarized by stating that the scheme of $n$-dimensional representations of a tree $\ST$ of semi-simple algebras decomposes into a finite disjoint union of its irreducible components
\[
\wis{rep}_n~\ST = \bigsqcup_{| \alpha | =n} \wis{rep}_{\alpha}~\ST = \bigsqcup_{| \alpha | = n} GL_n \times^{GL(\alpha)} \wis{rep}_{\alpha} U_{\alpha} \]
 where $\alpha$ is an $n$-dimension vector of $\Gamma(\ST)$ and where $U_{\alpha}$ is an affine Zariski open subset of the quiver-representation affine space $\wis{rep}_{\alpha}~\Gamma(\ST)$ determined by the determinant of the matrix obtained from the set of morphisms $\Sigma$. Moreover, the action of $GL_n$ encoding isomorphism classes is the same on both sides.
 
 In this section we will study the quotient varieties $\wis{iss}_n~\ST = \bigsqcup \wis{iss}_{\alpha}~\ST$ under the $GL_n$-action which classify the isomorphism classes of {\em semi-simple} $n$-dimensional representations of the tree $\ST$ of semi-simple algebras. To this end we introduce the {\em \'etale quiver} $\Psi = \Psi(\ST)$ see \cite{LBqaq}.
 
\begin{definition} Let $\ST$ be a tree of semi-simple algebras with Zariski quiver $\Gamma$ and corresponding Euler-form $\chi_{\Gamma}$. We construct the quiver $\Psi = \Psi(\ST)$ as follows
\begin{itemize}
\item{Vertices : let $\mathcal{S}$ be the commutative semigroup of all $n$-dimension vectors for $\Gamma$ (for all $n$) and let $\{ \alpha_1,\hdots,\alpha_k \}$ be the set of semigroup generators. Then, $\Psi$ will have $k$ vertices, vertex $i$ corresponding to the generator $\alpha_i$.}
\item{Arrows : the number of oriented arrows $\xymatrix{\vtx{i} \ar[r] & \vtx{j}}$ in $\Psi$ is given by the formula $\delta_{ij} - \chi_{\Gamma}(\alpha_i,\alpha_j)$.}
\end{itemize}
 We call $\Psi = \Psi(\ST)$ the {\em \'etale quiver} of the tree $\ST$ of semi-simple algebras.
\end{definition}

\begin{example} Continuing example~\ref{exgl2}, we calculate that the first component of the semigroup corresponding to $GL_2(\Z)$ has $8$ generators ($4$ for $n=1$ and $4$ for $n=2$)
\[
\begin{array}{c|c|cccc|ccc}
& n & a_1 & a_2 & a_3 & a_4 & b_1 & b_2 & b_6 \\
\hline
g_1 & 1 & 1 & 0 & 0 & 0 & 1 & 0 & 0 \\
g_2 & 2 & 1 & 1 & 0 & 0 & 0 & 0 & 1 \\
g_3 & 2 & 0 & 0 & 1 & 1 & 0 & 0 & 1 \\
g_4 & 2 & 1 & 0 & 0 & 1 & 0 & 0 & 1 \\
g_5 & 1 & 0 & 1 & 0 & 0 & 0 & 1 & 0 \\
g_6 & 2 & 0 & 1 & 1 & 0 & 0 & 0 & 1 \\
g_7 & 1 & 0 & 0 & 0 & 1 & 0 & 1 & 0 \\
g_8 & 1 & 0 & 0 & 1 & 0 & 1 & 0 & 0 \\
\end{array}
\]
and the second component gives an additional couple of generators (for $n=2$)
\[
\begin{array}{c|c|c|ccc}
& n & a_5 & b_3 & b_4 & b_5 \\
\hline
g_{9} & 2 & 1 & 1 & 1 & 0 \\
g_{10} & 2 & 1 & 0 & 0 & 1 \\
\end{array}
\]
Using the Euler-form of $\Gamma(GL_2(\Z))$ we deduce that the \'etale quiver $\Psi(GL_2(\Z))$ has two components, namely
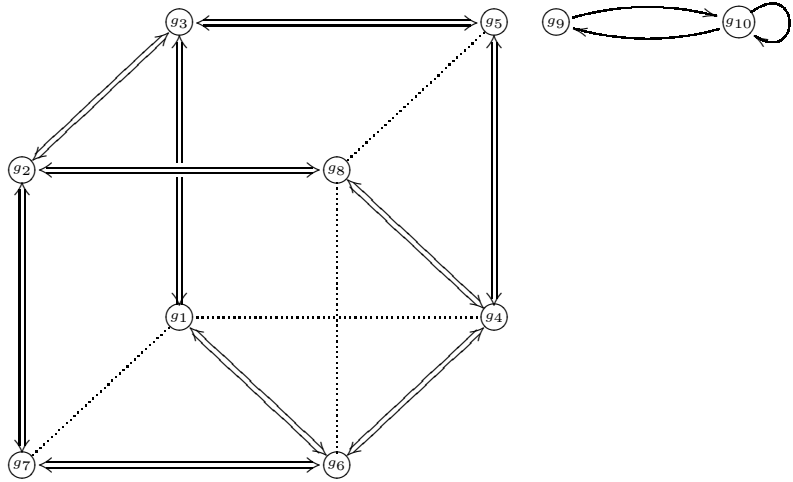
\[
\xymatrix@=.4cm{
& & & \vtx{g_3} \ar@{<=>}[rrrrrr] \ar@{<=>}'[ddd]'[dddddd] \ar@{<=>}[lllddd] & & & & & & \vtx{g_5} \\
& & & & & & & & & \\
& & & & & & & & & \\
\vtx{g_2} \ar@{<=>}[rrrrrr] \ar@{<=>}[dddddd] & & & & & & \vtx{g_8} \ar@{.}[rrruuu] \ar@{.}[dddddd] \ar@{<=>}[rrrddd] & & & \\
& & & & & & & & & \\
& & & & & & & & & \\
& & & \vtx{g_1} & & & & & & \vtx{g_4} \ar@{<=>}[uuuuuu] \ar@{.}[llllll] \\
& & & & & & & & & \\
& & & & & & & & & \\
\vtx{g_{7}} \ar@{<=>}[rrrrrr] \ar@{.}[rrruuu]  & & & & & & \vtx{g_6} \ar@{<=>}[llluuu]  \ar@{<=>}[rrruuu] & & & } \quad 
\xymatrix{ \vtx{g_9} \ar@/^/[rr] & & \vtx{g_{10}} \ar@/^/[ll] \ar@(ur,dr)}
\]
\end{example}

In these two examples we see that the \'etale quiver $\Psi$ is {\em symmetric}, that is, for all vertices $i,j$ we have
\[
\#~\{~\xymatrix{\vtx{i} \ar[r] & \vtx{j}}~\} = \#~\{~\xymatrix{\vtx{j} \ar[r] & \vtx{i}}~\}
\]
(though the number of {\em loops} in a vertex may be odd, as is the case in the second component of
$\psi(GL_2(\Z))$). This is not a coincidence.

\begin{theorem}\label{amalgamatedsym} For any tree $\ST$ of semi-simple algebras, the corresponding \'etale quiver $\Psi = \Psi(\ST)$ is symmetric.
\end{theorem}

\begin{proof}
We'll prove this theorem by induction on the number $z$ of edges of the tree. For $z=0$, the statement is trivial as the euler form of the quiver $\Gamma$ is $\chi_\Gamma=I$, and thus the quiver $\Psi$ also has euler form $\chi_\Psi=I$ so there are no arrows in $\Psi$.

Let the theorem be true for trees with $z-1$ edges. We now add a vertex $S'$ with edge $T$ to an existing vertex $S$,
$$\xymatrix{{\cdots} \ar@{-}[r] & \vtx{S} \ar@{-}[r]_{T} & \vtx{S'}}$$ 
with
$$S=\bigoplus_{i=1}^s \M_{l_i}(\C) \quad\text{ and }\quad
S'=\bigoplus_{i=1}^t \M_{m_i}(\C) \quad\text{ and }\quad
T=\bigoplus_{i=1}^r \M_{k_i}(\C).$$
Using all notations and conventions as before, we get
$$
\begin{bmatrix}l_1 \\ \vdots \\ l_s \end{bmatrix}
=
\underbrace{
\begin{bmatrix}
a_{11} & \hdots & a_{1r}\\
\vdots & \ddots & \vdots\\
a_{s1} & \hdots & a_{sr}
\end{bmatrix}
}_{A}
\begin{bmatrix}k_1 \\ \vdots \\ k_r\end{bmatrix}
\quad\text{ and }\quad
\begin{bmatrix}m_1 \\ \vdots \\ m_t \end{bmatrix}
=
\underbrace{
\begin{bmatrix}
b_{11} & \hdots & b_{1r}\\
\vdots & \ddots & \vdots\\
b_{t1} & \hdots & b_{tr}
\end{bmatrix}
}_{B}
\begin{bmatrix}k_1 \\ \vdots \\ k_r\end{bmatrix}.
$$
For a dimension vector 
$$\beta=(\hdots, \underbrace{p_1,\hdots, p_s}_{M|_S},\hdots,\underbrace{q_1,\hdots,q_t}_{M|_{S'}})$$
of an $n$-dimensional representation $M$ of the amalgamated product of the extended tree of semisimple algebras we get the following equality, using the restrictions $(M|_S)|_T$ and $(M|_{S'})|_T$,
\begin{equation}\label{dimveceq} \tag{$\dagger$}
\begin{bmatrix}p_1 & \cdots & p_s \end{bmatrix}
A
=
\begin{bmatrix}q_1 & \cdots & q_t \end{bmatrix}
B.
\end{equation}

This leads to the construction of the extended tree-partite quiver, giving a bipartite quiver for every edge of the tree, where the vertices corresponding to $S$ and $S'$ are their simple components and where the number of arrows $n_{ij}$ between the $i$'th vertex of $S$ and the $j$'th vertex of $S'$ is determined by the matrix
$$
D=
\begin{bmatrix}
n_{11} & \hdots & n_{1t}\\
\vdots & \ddots & \vdots\\
n_{s1} & \hdots & n_{st}
\end{bmatrix}
=
\begin{bmatrix}
a_{11} & \hdots & a_{1r}\\
\vdots & \ddots & \vdots\\
a_{s1} & \hdots & a_{sr}
\end{bmatrix}
\begin{bmatrix}
b_{11} & \hdots & b_{t1}\\
\vdots & \ddots & \vdots\\
b_{1r} & \hdots & b_{tr}
\end{bmatrix}
=A\ B\tp
$$
Adding the vertex $S'$ extends the euler form $\chi$ of the tree-partite quiver on $z-1$ edges downwards and to the right, giving the euler form $\chi'$ for the tree-partite quiver for the extended tree
$$
\chi'=
\left[
\begin{tabular}{C C C | C}
& & & 0 \\
& & & \vdots \\
& & & 0 \\
& \chi & & -D \\
& & & 0 \\
& & & \vdots \\
& & & 0 \\
\hline
0 &  \cdots & 0 & I_t
\end{tabular}
\right]
$$
where the row of the matrix $D$ corresponds to the existing vertex $S$.

We now calculate the number of arrows in the one quiver $\Psi$ corresponding to this euler form $\chi'$ between dimension vectors
$$\beta=(\hdots,p_1,\hdots, p_s,\hdots,q_1,\hdots,q_t)$$
and
$$\beta'=(\hdots,p'_1,\hdots, p'_s,\hdots,q'_1,\hdots,q'_t)$$
using $\chi'(\beta,\beta')=\beta\cdot\chi'\cdot\beta'$,
$$\beta\cdot\chi'\cdot\beta'
= \beta_{\setminus q}\cdot\chi\cdot\beta'_{\setminus q'} 
\ -\ 
\begin{bmatrix}p_1 & \cdots & p_s\end{bmatrix}
D
\begin{bmatrix}q'_1 \\ \vdots \\ q'_t\end{bmatrix}
\ +\ 
\begin{bmatrix}q_1 & \cdots & q_t\end{bmatrix}
\begin{bmatrix}q'_1 \\ \vdots \\ q'_t\end{bmatrix}.
$$
The difference between the number of arrows in one direction and arrows in the other direction should be zero for the quiver $\Psi$ to be symmetric. Indeed,
$$
\beta\cdot\chi'\cdot\beta'-\beta'\cdot\chi'\cdot\beta
=
-
\begin{bmatrix}p_1 & \cdots & p_s\end{bmatrix}
D
\begin{bmatrix}q'_1 \\ \vdots \\ q'_t\end{bmatrix}
+
\begin{bmatrix}p'_1 & \cdots & p'_s\end{bmatrix}
D
\begin{bmatrix}q_1 \\ \vdots \\ q_t\end{bmatrix}
$$
is the sum of two scalars so taking the transpose of the second summand makes this equal to
$$
-
\begin{bmatrix}p_1 & \cdots & p_s\end{bmatrix}
D
\begin{bmatrix}q'_1 \\ \vdots \\ q'_t\end{bmatrix}
+
\left(
\begin{bmatrix}p'_1 & \cdots & p'_s\end{bmatrix}
D
\begin{bmatrix}q_1 \\ \vdots \\ q_t\end{bmatrix}
\right)\tp
$$
or
$$
-
\begin{bmatrix}p_1 & \cdots & p_s\end{bmatrix}
A\ B\tp
\begin{bmatrix}q'_1 \\ \vdots \\ q'_t\end{bmatrix}
+
\begin{bmatrix}q_1 & \cdots & q_t\end{bmatrix}
B\ A\tp
\begin{bmatrix}p'_1 \\ \vdots \\ p'_s\end{bmatrix}
$$
which is $0$ according to equation (\ref{dimveceq}).
\end{proof}

\begin{remark} \label{rem2}
The results remains true if we allow vertex algebras $S_v \simeq \C[x]$ or $\C[x,x^{-1}]$. The Euler form of the relevant part is then
\[
\chi = \begin{bmatrix} 0 & -1 & \hdots & -1 \\
0 & & & \\
\vdots & & I_k & \\
0 & & & \end{bmatrix} \]
and the dimension vectors $\alpha = (a_0,a_1,\hdots,a_k)$ and $\beta=(b_0,b_1,\hdots,b_k)$ satisfy $a_0 = \sum_{i=1}^k a_i$ and $b_0=\sum_{i=1}^k b_i$ whence
\[
\chi(\alpha,\beta) = \sum_{i=1}^k a_ib_i -a_0b_0 = \chi(\beta,\alpha) \]
\end{remark}

A first application of the \'etale quiver $\Psi(\ST)$ is that allows us to determine the components $\wis{rep}_{\alpha}~\ST$ containing a (Zariski open subset of) simple representation(s) of $\ST$. Recall that in \cite{LBProcesi} a characterization was proved of the dimension vectors of simple representations of quivers. Remark that if $\alpha_i$ is a generator of the semigroup $\mathcal{S}$ of non-empty components in the representation varieties of $\ST$, then all points of $\wis{rep}_{\alpha_i}~\ST$ are simple representations. As a consequence, for any component $\wis{rep}_{\alpha}~\ST$ we can write
\[
\alpha = n_1 \alpha_1 + \hdots + n_k \alpha_k \]
with all $n_i \in \N$ and hence $\wis{rep}_{\alpha}~\ST$ contains semi-simple representations
\[
M = S_1^{\oplus n_1} \oplus \hdots \oplus S_k^{\oplus n_k} \]
with $S_i$ a simple representation in $\wis{rep}_{\alpha_i}~\ST$.

\begin{theorem} Let $\ST$ be a tree of semi-simple algebras with associated \'etale quiver $\Psi$ and Euler-form $\chi_{\Psi}$. Let $\{ \alpha_1,\hdots,\alpha_k \}$ be the set of generators of the semigroup of non-empty components $\wis{rep}_{\alpha}~\ST$ of the representation varieties of $\ST$. Then, the following statements are equivalent
\begin{enumerate}
\item{The component $\wis{rep}_{\alpha}~\ST$ contains a simple representation.}
\item{The component $\wis{rep}_{\alpha}~\ST$ contains a non-empty Zariski open subset of simple representations.}
\item{For a decomposition $\alpha = n_1 \alpha_1 + \hdots + n_k \alpha_k$ with all $n_i \in \N$ we have that $\beta = (n_1,\hdots,n_k)$ is the dimension vector of a simple representation of $\Psi$.}
\item{For any decomposition $\alpha = n_1 \alpha_1 + \hdots + n_k \alpha_k$ with all $n_i \in \N$ we have that $\beta = (n_1,\hdots,n_k)$ is the dimension vector of a simple representation of $\Psi$.}
\end{enumerate}
Recall from \cite{LBProcesi} that $\beta$ is the dimension vector of a simple representation of $\Psi$ if and only if $\wis{supp}(\beta)$ is a strongly connected subquiver of $\Psi$ and for all vertices $i \in \wis{supp}(\beta)$ we have $\chi_{\Psi}(\epsilon_i,\beta) \leq 0$ and $\chi_{\Psi}(\beta,\epsilon_i) \leq 0$ unless $\wis{supp}(\beta) = \{ i,j \}$ or $\{ i \}$ and there is just one arrow from $i$ to $j$ (and one form $j$ to $i$ by symmetry) (resp. one loop in $i$) in which case the non-zero components of $\beta$ are equal to one.
\end{theorem}

\begin{proof} Compare with \cite[Prop. 6]{LBqaq}. Equivalences $(1)$ and $(2)$ follow from Schur's lemma. Equivalences $(3)$ and $(4)$ follow from the fact that the component $\wis{rep}_{\alpha}~\ST$ is an irreducible variety.
\end{proof}

A point $\xi$ of the quotient variety $\wis{iss}_{\alpha}~\ST = \wis{rep}_{\alpha}~\ST / GL_n$ determines a semi-simple representation
\[
M_{\xi} = T_1^{\oplus e_1} \oplus \hdots \oplus T_l^{\oplus e_l} \]
where each $T_i \in \wis{rep}_{\beta_i}~\ST$ is a simple $\ST$-representation. Observe that by the foregoing result we have full information on the dimension vectors $\beta_i$ which can occur in such a decomposition. Each of the $\beta_i$ can be written (possibly in several different ways) as a linear combination of the generators $\alpha_j$
\[
\beta_i = c^{(i)}_1 \alpha_1 + \hdots + c^{(i)}_k \alpha_k \qquad \text{with $c^{(i)}_j \in \N$.} \]
Define the dimension vectors (of the \'etale quiver $\Psi = \Psi(\ST)$) $\gamma_i = (c^{(i)}_1,\hdots,c^{(i)}_k)$.

\begin{definition} The {\em local quiver setting} of $M_{\xi}$ is the quiver $\Delta_{\xi}$ and the dimension vector $\alpha_{\xi}$ where $\Delta_{\xi}$ has
\begin{itemize}
\item{$l$ vertices $\{ v_1,\hdots,v_l \}$ corresponding to the non-isomorphic simple components of $M_{\xi}$, and}
\item{$\delta_{ij} - \chi_{\Psi}(\gamma_i,\gamma_j)$ directed arrows $\xymatrix{\vtx{v_i} \ar[r] & \vtx{v_j}}$.}
\item{The dimension vector $\alpha_{\xi} = (e_1,\hdots,e_l)$ corresponding to the multiplicities of the simple components in $M_{\xi}$.}
\end{itemize}
Remark that as $\Psi = \Psi(\ST)$ is symmetric, so is the quiver $\Delta_{\xi}$. Observe that all information about the possible local quiver settings is contained in the \'etale quiver $\Psi(\ST)$.
\end{definition}

The main use for the local quiver settings is that they determine the \'etale local structure of the quotient variety $\wis{iss}_{\alpha}~\ST$ in the point $\xi$ (which also explains the terminology {\em \'etale quiver} for $\Psi$).

\begin{theorem} \label{etale} There is an \'etale isomorphism between
\begin{itemize}
\item{$\wis{iss}_{\alpha}~\ST$ in the point $\xi$ and}
\item{$\wis{iss}_{\alpha_{\xi}}~\Delta_{\xi}$ in the point $\overline{0}$ determined by the zero representation.}
\end{itemize}
\end{theorem}

\begin{proof} See \cite[Thm 4]{LBqaq}.
\end{proof}

In particular, this result can be used to determine the smooth locus of the quotient varieties $\wis{iss}_{\alpha}~\ST$ using the fact all local quivers are symmetric and the characterization, due to Raf Bocklandt \cite{Bocklandt1}, of all coregular symmetric quiver settings.

\begin{definition} Let $\Delta$ be a symmetric quiver and $\alpha$ a dimension vector. We say that the quiver setting $(\Delta,\alpha)$ is a {\em cherry-tree} if the following conditions are satisfied
\begin{itemize}
\item{Replace in $\Delta | \wis{supp}(\alpha)$ each symmetric-arrow pair by one solid line, then (multiple) cherries can grow only at a node where $\alpha$ is equal to one and where one cherry looks like
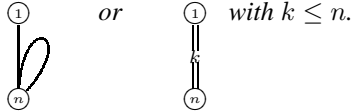
\[
\xymatrix{\vtx{1} \ar@{-}[d] \\ \vtx{n} \ar@{-}@(ur,u)}~\qquad \text{ or}~\qquad
\xymatrix{\vtx{1} \ar@{=}[d]|{k} \\ \vtx{n}}~\quad~\text{ with $k \leq n$.}
\]
}
\item{After harvesting all cherries, we are left with a tree made from mini-branches of  the following forms
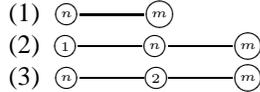
\begin{enumerate}
\item{$\xymatrix{\vtx{n} \ar@{-}[r] & \vtx{m}}$}
\item{$\xymatrix{\vtx{1} \ar@{-}[r] & \vtx{n} \ar@{-}[r] & \vtx{m}}$}
\item{$\xymatrix{\vtx{n} \ar@{-}[r] & \vtx{2} \ar@{-}[r] & \vtx{m}}$}
\end{enumerate}
joined at common nodes where $\alpha$ is equal to one.}
\end{itemize}
\end{definition}

\begin{theorem}[Bocklandt,\cite{Bocklandt1}] The following are equivalent
\begin{enumerate}
\item{The quotient variety $\wis{iss}_{\alpha}~\Delta$ is smooth.}
\item{The quiver setting $(\Delta,\alpha)$ is a cherry-tree.}
\end{enumerate}
\end{theorem}

As all local quivers $\Delta_{\xi}$ are symmetric, theorem~\ref{etale} and Bocklandt's result allow us to verify whether $\xi \in \wis{iss}_{\alpha}~\ST$ is a smooth point. Further, we can use the \'etale quiver $\Psi = \Psi(\ST)$ to determine the components $\wis{iss}_{\alpha}~\ST$ which are smooth affine varieties.

\begin{theorem} The following are equivalent
\begin{enumerate}
\item{The quotient variety $\wis{iss}_{\alpha}~\ST$ is smooth.}
\item{For any decomposition $\alpha = n_1 \alpha_1+\hdots+n_k \alpha_k$ with $n_i \in \N$ and $\alpha_i$ the generators of $\mathcal{S}$, we have that $(\Psi,\beta)$ is a cherry-tree where $\beta = (n_1,\hdots,n_k)$.}
\end{enumerate}
\end{theorem}

\begin{proof} For a decomposition $\alpha = n_1 \alpha_1+\hdots+n_k \alpha_k$ we obtain points $\xi$ in $\wis{iss}_{\alpha}~\ST$ with corresponding semi-simple module
\[
M_{\xi} = S_1^{\oplus n_1} \oplus \hdots \oplus S_k^{\oplus n_k} \qquad \text{with $S_i$ a simple $\alpha_i$-dimensional representation} \]
and one verifies that the local quiver setting in such a point $(\Delta_{\xi},\alpha_{\xi}) = (\Psi | \wis{supp}(\beta), \beta)$. Hence, $(1)$ implies $(2)$. Conversely, $(2)$ implies that $\wis{iss}_{\alpha}~\ST$ is smooth in the points $\xi$ constructed as above. As any other point degenerates to one of the $\xi$'s, see for example \cite{LBProcesi}, the result follows. Alternatively, if $\psi$ is another point, its local quiver setting can be obtained from one of the $(\Delta_{\xi},\alpha_{\xi})$ above and smoothness is preserved under taking local quivers.
\end{proof}

\section{Double Poisson structure}
 
 In the foregoing section we have seen that the local structure of the quotient varieties $\wis{iss}_{\alpha}~\ST$ near a point $\xi$ is determined by a symmetric local quiver setting $\wis{iss}_{\alpha_{\xi}}~\Delta_{\xi}$. In this section we will prove that these quiver quotient varieties have a lot of additional structure, in particular, they carry a natural Poisson structure induced by a double Poisson algebra structure, as introduced by Michel Van den Bergh in \cite{VdBPoisson}, on a suitable universal localization of the path algebra of a quiver related to the \'etale quiver.
 
\begin{definition} From \cite{VdBPoisson} we recall that a $\C$-algebra $A$ is said to be a {\em double Poisson algebra} if there is a bilinear map
\[
\ldb -,- \rdb~:~A \times A \rTo A \otimes A \]
 satisfying the following properties for all $a,b,c \in A$
\begin{enumerate}
\item{$\ldb a,bc \rdb = \ldb a,b \rdb c + b \ldb a,c \rdb$.}
\item{$\ldb a,b \rdb = - \ldb b,a \rdb^o$ where $(u \otimes v)^o = v \otimes u$.}
\item{$\ldb a,b,c \rdb = \ldb a, \ldb b,c \rdb \rdb_L + \tau_{(123)} \ldb b , \ldb c,a \rdb \rdb_L + \tau_{(132)} \ldb c,\ldb a,b \rdb \rdb_L = 0$ where if $\ldb b,c \rdb = u_{(1)} \otimes u_{(2)}$ then $\ldb a,\ldb b,c \rdb \rdb_L = \ldb a,u_{(1)} \rdb \otimes u_{(2)}$.}
\end{enumerate}
\end{definition}

If $A$ is a double Poisson algebra, there are induced Poisson structures on all its representation schemes $\wis{rep}_n~A$ and corresponding quotient schemes $\wis{iss}_n~A$. Observe that any $a \in A$ defines matrix-valued functions $a_{ij} \in \C[\wis{rep}_n~A]$ and by \cite[Prop 1.2]{VdBPoisson} the Poisson bracket on $\C[\wis{rep}_n~A]$ is defined by
\[
\{ a_{ij},b_{uv} \} = (\ldb a,b \rdb_{(1)})_{uj} (\ldb a,b \rdb_{(2)})_{iv} \]
Recall from \cite{LBProcesi} that $\C[\wis{iss}_n~A]$ is generated by the functions $a_{ii}=tr(a)$ for $a \in A$ and from \cite[Prop 7.7.2]{VdBPoisson} we recall that there is a unique induced Poisson structure on the quotient variety $\wis{iss}_n~A$ coming from the Poisson structure on $\C[\wis{iss}_n~A]$ given by
\[
\{ tr(a),tr(b) \} = tr \{ a,b \} \]
Observe that these brackets define plenty of additional structure (differential forms, poly vectorfields, Schouten brackets, Hamiltonian structure etc.) on the representation schemes $\wis{rep}_n~A$ and their quotient schemes $\wis{iss}_n~A$, see \cite[\S 7]{VdBPoisson} for more details.

In this paper, we are only interested in the canonical double Poisson structure defined by Michel Van den Bergh on the path algebra $A = \C \Delta^d$ of a {\em double quiver} $\Delta^d$.

\begin{definition} If $\Delta$ is a finite quiver on vertices $\{ 1,\hdots,l \}$, then its {\em double quiver} $\Delta^d$ is the extended quiver obtained from $\Delta$ by adjoining for every arrow $a$ (resp. loop $l$) in $\Delta$ a new arrow $a^*$ (resp. loop $l^*$) in the reverse direction, that is, if
\[
\xymatrix{\vtx{i} \ar[r]^a & \vtx{j}} \qquad \text{then we have } \qquad \xymatrix{\vtx{i} & \vtx{j} \ar[l]^{a^*}} \]
and for any loop $l$ in $\Delta$ 
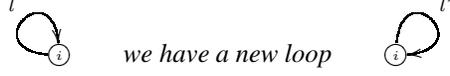
\[
\xymatrix{\vtx{i} \ar@(l,u)^l} \qquad \text{we have a new loop} \qquad \xymatrix{\vtx{i} \ar@(u,r)^{l^*}} \]
Observe that a symmetric quiver is a double quiver if only if the number of loops in every vertex is even.
\end{definition}

Let $V = \overbrace{\C \times \hdots \times \C}^l$ be the subalgebra of $\C \Delta^d$ generated by the vertex-idempotents $e_i$. For every arrow (or loop) $\xymatrix{\vtx{i} \ar[r]^c & \vtx{j}}$ in $\Delta^d$ we define the $V$-derivation
\[
\frac{\partial}{\partial c} \in Der_V(\C \Delta^d,\C \Delta^d \otimes \C \Delta^d) \quad \text{by} \quad \frac{\partial b}{\partial c} = \begin{cases} e_j \otimes e_i & \text{if $b=c$} \\ 0 & \text{otherwise} \end{cases}
\]
Consider the canonical differential $2$-form (the Hamiltonian structure of \cite[\S 6.3]{VdBPoisson})
\[
P = \sum_{a \in \Delta} \frac{\partial}{\partial a} \frac{\partial}{\partial a^*} \]
then this defines a double Poisson algebra structure on the path algebra $\C \Delta^d$ by taking for all paths $p,q \in \C \Delta^d$
\[
\ldb p,q \rdb = \sum_{a \in \Delta} \frac{\partial q}{\partial a^*}_{(1)} \frac{\partial p}{\partial a}_{(2)} \otimes
\frac{\partial p}{\partial a}_{(1)} \frac{\partial q}{\partial a^*}_{(2)} - \frac{\partial q}{\partial a}_{(1)} \frac{\partial p}{\partial a^*}_{(2)} \otimes \frac{\partial p}{\partial a^*}_{(1)} \frac{\partial q}{\partial a}_{(2)} \]
This double Poisson structure then induces canonical Poisson structures on all representation varieties $\wis{rep}_{\alpha}~\Delta^d$ and their quotient varieties $\wis{iss}_{\alpha}~\Delta^d$.

For a tree $\ST$ of semi-simple algebras, the local structure of the quotient varieties $\wis{iss}_{\alpha}~\ST$ is given by a symmetric (rather than a double) quiver. Still, this is only a minor annoyance.

\begin{theorem} For any symmetric quiver setting $(\Delta,\alpha)$ there is a canonical Poisson structure on the quotient variety
\[
\wis{iss}_{\alpha}~\Delta \]
induced by a double Poisson algebra structure on a universal localization of the path algebra $\C \tilde{\Delta}$ of an associated double quiver $\tilde{\Delta}$.
\end{theorem}

\begin{proof} As $\Delta$ is symmetric all non-loop arrows appear in pairs $\{ a,a^* \}$ so we only have to worry about unaccompanied loops $l$ in certain vertices $i$. The trick is to extend the quiver $\Delta$ for each such loop by adding a new vertex $i'$ and replacing the loop $l$ by a pair of arrows $\{ a_l,a_l^* \}$
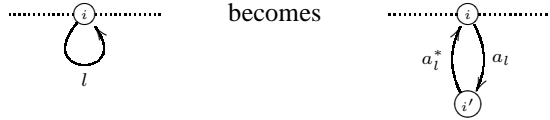
\[
\xymatrix{\ar@{.}[r] & \vtx{i} \ar@(dl,dr)_l \ar@{.}[r] &} \qquad \text{becomes} \qquad 
\xymatrix{
\ar@{.}[r] & \vtx{i} \ar@/^/[d]^{a_l} \ar@{.}[r] & \\
& \vtx{i'} \ar@/^/[u]^{a_l^*} & \\}
\]
If we repeat this procedure for all loose loops we end up with a double quiver $\tilde{\Delta}$. By the above remarks, its path algebra $\C \tilde{\Delta}$ carries a canonical double Poisson algebra structure.
In $\C \tilde{\Delta}$ consider the set of arrows $\Sigma = \{ a_l~:~l~\text{a loose loop}~\}$. 

By \cite[Prop 2.4.3]{VdBPoisson} the double Poisson bracket $\ldb -,- \rdb$ on $\C \tilde{\Delta}$ extends uniquely to a double Poisson bracket on the universal localization $A = \C \tilde{\Delta}_{\Sigma}$. As the matrices corresponding to the arrows $a_l$ in a finite dimensional representation of $A$ must be invertible it follows that $\wis{rep}_{\beta}~A$ is empty unless $\beta = \tilde{\alpha}$ whenever $\tilde{\alpha}(i') = \tilde{\alpha}(i)$. Moreover, the map
\[
\wis{rep}_{\tilde{\alpha}}~A \rTo \wis{rep}_{\alpha}~\Delta \]
sending all matrix-couples $(M(a_l),M(a_l^*))$ to $M(a_l^*)M(a_l)$ is a principal fibration with group the basechange group in the additional vertices. As a consequence, we have that
\[
\wis{iss}_{\tilde{\alpha}}~A \simeq \wis{iss}_{\alpha}~\Delta \]
finishing the proof.
\end{proof}

\begin{example} The canonical \'etale Poisson structure on the quotient varieties $\wis{iss}_{\alpha}~GL_2(\Z)$ is induced from the double Poisson algebra structure on the double quiver $\Delta^d$ where $\Delta$ is the quiver
\[
\xymatrix@=.4cm{
& & & \vtx{g_3} \ar[rrrrrr] \ar@{<-}'[ddd]'[dddddd] \ar@{<-}[lllddd] & & & & & & \vtx{g_5} \\
& & & & & & & & & \\
& & & & & & & & & \\
\vtx{g_2} \ar@{<-}[rrrrrr] \ar@{<-}[dddddd] & & & & & & \vtx{g_8} \ar@{.}[rrruuu] \ar@{.}[dddddd] \ar@{<-}[rrrddd] & & & \\
& & & & & & & & & \\
& & & & & & & & & \\
& & & \vtx{g_1} & & & & & & \vtx{g_4} \ar@{<-}[uuuuuu] \ar@{.}[llllll] \\
& & & & & & & & & \\
& & & & & & & & & \\
\vtx{g_{7}} \ar@{<-}[rrrrrr] \ar@{.}[rrruuu]  & & & & & & \vtx{g_6} \ar@{->}[llluuu]  \ar@{<-}[rrruuu] & & & } \quad 
\xymatrix{ \vtx{g_9} \ar[rr] & & \vtx{g_{10}} \ar[dd]|{\bullet} \\
& & \\
& & \vtx{g_{11}}}
\]
and extended to the universal localization by inverting the $\bullet$-arrow.
\end{example}

If $\ldb -,- \rdb$ is a double Poisson bracket on $A$ we recall from \cite[Prop 1.4]{VdBPoisson} that the induced bracket
\[
\{ -,- \}~:~A \times A \rTo A \qquad (a,b) \mapsto \ldb a,b \rdb_{(1)} \ldb a,b \rdb_{(2)} \]
is a derivation in the second argument and vanishes on commutators in the first argument. Moreover, this bracket turns $A$ into a left {\em Loday algebra} meaning that the following relation holds for all $a,b,c \in A$
\[
\{ a, \{ b,c \} \} = \{ \{ a,b \},c \} + \{ b, \{ a,c \} \} \]
This, in turn, implies that the vectorspace quotient $\wis{neck} = A/[A,A]$ is a Lie algebra, the so called {\em necklace Lie algebra} introduced in \cite{LBBocklandt} and \cite{Ginzburg}.

\section{Necklaces and flows}
 
 In the foregoing section we proved that there is an \'etale Poisson structure on the quotient varieties $\wis{iss}_{\alpha}~\ST$ and that the double Poisson algebra defines an infinite dimensional Lie algebra $\wis{neck}$, the {\em necklace Lie algebra}. In this section we will see that these necklaces induce vectorfields and flows on the quiver quotient varieties, thereby inducing a dynamic aspect on all $\wis{iss}_{\alpha}~\ST$.

As seen before (theorem~\ref{amalgamatedsym}), a large class of algebras has symmetric (potentially double) quivers $\Delta$ associated to them. The involution $*$ gives rise to a symplectic structure on the path algebra $\C{\Delta}$.. The study of non-commutative differential forms for path algebras of quivers was independently developed in \cite{LBBocklandt} and \cite{Ginzburg}. A concise outline of this construction is given here, following the first article.

\begin{definition}
Let $\Delta$ be a quiver on $k$ vertices. A basis for the relative
differential $n$-forms $\Omega_v^n \Delta$ is given by the elements $$p_0\ dp_1\ \cdots\ dp_n$$ where the $p_i$ are oriented paths in $\Delta$, which are nontrivial for $i\geq 1$, and such that the starting point of $p_i$ is the end point of $p_{i+1}$. 
\end{definition}

\begin{definition}
A \emph{necklace} or \emph{necklace word} $n$ in $\Delta$ is an oriented cyclic path in $\Delta$ upto cyclic equivalence.

If $a$ is an arrow in a necklace $n$, we can define the \emph{derivative $\frac{\partial n}{\partial a}$} as the sum of all oriented paths
obtained by deleting one occurrence of $a$ in $n$.
The \emph{differential $d$} takes a necklace $n$ in $\Delta$ to
$$dn\defeq \sum_{a\in\arr}\frac{\partial n}{\partial a}da.$$
\end{definition}

After dividing out the supercommutators,
$$\dR_v^n \Delta \defeq \frac{\Omega_v^n \Delta}{\sum_{i=1}^n \left[\Omega_v^i \Delta,\ \Omega_v^{n-i} \Delta\right]},$$
we obtain the relative Karoubi complex,
$$\dR_v^0 \overset{d}{\map} \dR_v^1 \Delta \overset{d}{\map} \dR_v^2 \Delta \overset{d}{\map} \cdots$$
The first of these terms $\dR_v^i$ are known explicitly.
$\dR_v^0 \Delta=\frac{\C \Delta}{[\C \Delta, \C \Delta]}$ is the vectorspace spanned by all \emph{necklace words} in $\Delta$.
$\dR_v^1 \Delta$ is the vectorspace defined by $$\dR_v^1=\setd{\bigoplus_{a\in\arr}p\ da}{p \text{ is a path in } \Delta \text{ starting in } t(a)}.$$

Let $\Delta$ now be a \emph{double quiver} and let $L\subset \arr$ be the equivalence class of all unstarred arrows in $\arr$.

\begin{definition}
For each double quiver $\Delta$, we can define a Lie algebra structure on the space $\dR_v^0 A$ by
$$[n,n']\defeq \sum_{a\in L}\left( \frac{\partial n}{\partial a}\ \frac{\partial n'}{\partial a^*} - \frac{\partial n}{\partial a^*}\ \frac{\partial n'}{\partial a}\right).$$
We call this the \emph{necklace Lie algebra} of $\Delta$.
\end{definition}

This can be illustrated by a simplified example sketch. Let $n$ and $n'$ be two necklaces having where $n$ contains exactly one un-$*$-ed arrow $a$ for which $n'$ contains the $*$-ed $a^*$, and $n$ does not contain $a^*$ for which $n'$ contains $a$, then 
$$
[n,n']=
\left[
\centerprent{\xymatrix@=1em{
\\
& \vtx{} \ar@/^.2em/[rr]^{a} & & \vtx{} \ar@/^.2em/[rd] \\
\vtx{} \ar@/^.2em/[ru] & & & & \vtx{} \ar@{-->}@/^/[llll]
}},
\centerprent{\xymatrix@=1em{
\vtx{} \ar@{-->}@/^/[rrrr] & & & & \vtx{} \ar@/^.2em/[ld] \\
& \vtx{} \ar@/^.2em/[lu] & & \vtx{} \ar@/^.2em/[ll]^{a^*} \\ {}
}}
\right]
=
\centerprent{\xymatrix@=1em{
\vtx{} \ar@{-->}@/^/[rrrr] & & & & \vtx{} \ar@/^.2em/[ld] \\
& \vtx{} \ar@/^.2em/[lu] & & \vtx{} \ar@/^.2em/[rd] \\
\vtx{} \ar@/^.2em/[ru] & & & & \vtx{} \ar@{-->}@/^/[llll]
}} -\quad 0
$$
and for multiple occurrences of couples of $*$-ed and un-$*$-ed arrows in these necklaces, sums of the above operations are to be taken. Examples of these Lie brackets can be found in the examples~\ref{ex:liepsl2} and \ref{ex:liegl2}.

Let $\dR_v^1 \Delta_{ex}$ be the kernel of $\dR_v^1 \Delta \overset{d}{\map} \dR_v^2 \Delta$.

\begin{definition}
We define $\Der_\omega \Delta$ the Lie algebra of \emph{symplectic derivations} of $\C \Delta$ as the $\C^{\times k}$-derivations of $\C \Delta$ which preserve the moment element $$m\defeq \sum_{a\in L}aa^*-a^*a\in\C \Delta,$$
being the sum of all commutators of couples of $*$-ed and corresponding un-$*$-ed arrows.
\end{definition}

\begin{theorem}[{\cite[th.~4.2]{LBBocklandt}}]
We can now identify $\dR_v^1 \Delta_{ex}$ with this $\Der_\omega \Delta$ by assigning to the image $dn\in\dR_v^1 \Delta_{ex}$ of a necklace $n$ the $\C^{\times k}$-derivation $\theta_n\in\Der_\omega \Delta$, defined by
$$\theta_n(a)\defeq-\frac{\partial n}{\partial a^*}\qquadtext{ and }\theta_n(a^*)\defeq\frac{\partial n}{\partial a}.$$
\end{theorem}

For example, this symplectic structure can be put on all \emph{double} quivers arising out of theorem~\ref{amalgamatedsym}, e.g.\ those of examples~\ref{ex:liepsl2} and \ref{ex:liegl2}, i.e.\
$$
\centerprent{\xymatrix@=1em{
& \vtx{} \ar@/^.2em/[dl]^{} \ar@/^.2em/[dr]^{} & \\
\vtx{}  \ar@/^.2em/[ur]^{} \ar@/^.2em/[d]^{} & & \vtx{}  \ar@/^.2em/[ul]^{} \ar@/^.2em/[d]^{} \\
\vtx{}  \ar@/^.2em/[dr]^{} \ar@/^.2em/[u]^{} & & \vtx{}  \ar@/^.2em/[dl]^{} \ar@/^.2em/[u]^{} \\
& \vtx{} \ar@/^.2em/[ur]^{} \ar@/^.2em/[ul]^{} 
}}
\qquadtext{ and } 
\centerprent{
\xymatrix@=1.5em{
\vtx{~} \ar@/^.2em/[r]^{} \ar@/^.2em/[d]^{} &
\vtx{~} \ar@/^.2em/[l]^{} \ar@/^.2em/[r]^{} \ar@/^.2em/[dd]^(.75){} &
\vtx{~} \ar@/^.2em/[l]^{} \ar@/^.2em/[d]^{} \\
\vtx{~} \ar@/^.2em/[u]^{} \ar@/^.2em/[rr]^(.75){} \ar@/^.2em/[d]^{} & &
\vtx{~} \ar@/^.2em/[ll]^(.75){} \ar@/^.2em/[u]^{} \ar@/^.2em/[d]^{} \\
\vtx{~} \ar@/^.2em/[u]^{} \ar@/^.2em/[r]^{} &
\vtx{~} \ar@/^.2em/[uu]^(.75){} \ar@/^.2em/[l]^{} \ar@/^.2em/[r]^{} &
\vtx{~} \ar@/^.2em/[l]^{} \ar@/^.2em/[u]^{}
}}.
$$

\begin{definition}
A derivation $\theta_n \in \Der_\omega \Delta$ is \emph{locally nilpotent} if for any path $p \in \C{\Delta}$ there exists a $k\in \N$ such that $\theta_n^k(p)=0$.
\end{definition}

\begin{definition}
A necklace $n$ is called a \emph{one-way necklace} if it does not contain both $a$ and $a^*$, i.e.\ if for all arrows $a$ we have $\set{a,a^*}\not\subseteq\Supp(n)$.
\end{definition}

\begin{lemma}\label{le_locallynilpotent}
A derivation $\theta_n$ is locally nilpotent if and only if $n$ is a one-way necklace.
\end{lemma}

\begin{proof}
{[$\pmi$]}. Let $p$ be a path in $\C{\Delta}$ and suppose $p=a_l\cdots a_1$ contains $k$ arrows not in $\Supp(n)$. Now every term of $\theta_n(p)$,
$$a_l\cdots \theta_n(a_i)\cdots a_1$$
contains $k-1$ arrows not in $\Supp(n)$, because 
\begin{itemize}
\item if $a_i\in\Supp(n)$, then $\theta_n(a_i)=0$ because $a_i^*\nin\Supp(n)$;
\item if $a_i\nin\Supp(n)$, then $\Supp(\theta_n(a_i))\subseteq \Supp(n)$.
\end{itemize}
This yields $$\theta_n^{k+1}(p)=0.$$

{[$\imp$]}. Let $n$ be a necklace containing $k$ times the arrow $a$ and $l$ times the arrow $a^*$,  with $k,l\geq 1$. We denote this degree by $\#_{a}(n)=k$ and $\#_{a^*}(n)=l$. So all $l$ terms $p$ of $\theta_n(a)$ are negative and 
$$\#_{a}(p)=k \qquadtext{ and } \#_{a^*}(p)=l-1.$$
All $k$ terms $p'$ of $\theta_n(a^*)$ are positive and
$$\#_{a}(p')=k-1 \qquadtext{ and } \#_{a^*}(p')=l.$$
For every term $p''$ of $\theta_n(z)$ for any other arrow $z\in\Supp(n)$ we have
$$\#_{a}(p'')=k \qquadtext{ and } \#_{a^*}(p'')=l.$$
We will now study the terms of the lowest degree in $a$ and $a^*$ of $\theta_n^t(a^*)$ for $t\in \N$. 

\begin{itemize}
\item
Considering $\theta_n^2(a^*)$, the terms of the lowest degree are those of the form
$$z_x\cdots \theta_n(a^*)\cdots z_1\qquadtext{ and }z_x\cdots\theta_n(a)\cdots z_1$$
with $z_x\cdots z_1$ being a term (of minimal degree in $a$ and $a^*$) of $\theta_n(a^*)$. Both terms contain $2k-2$ times the arrow $a$ and $2l-1$ times the arrow $a^*$. Moreover, the first term has positive sign, the second negative. There are $k^2l$ terms of the first kind ($k$ terms $z_x\cdots z_1$, $k$ times $a$ in $n$ and $l$ times $a^*$ in $z_x\cdots z_1$), and $kl(k-1)$ of the second kind ($k$ terms $z_x\cdots z_1$, $l$ times $a^*$ in $n$ and $k-1$ times $a$ in $z_x\cdots z_1$). So if the derivation is locally nilpotent of degree two for $a^*$, then these kinds of terms need to cancel each other out, so $k^2l-kl(k-1)=kl=0$, which is impossible. 

\item
Terms of the lowest degree of $\theta_n^3(a^*)$ are again those of the form 
$$z_x\cdots \theta_n(a^*)\cdots z_1\qquadtext{ and }z_x\cdots\theta_n(a)\cdots z_1$$
with $z_x\cdots z_1$ being a term of minimal degree in $a$ and $a^*$ of $\theta_n^2(a^*)$. Both above terms contain $3k-3$ times the arrow $a$ and $3l-2$ times the arrow $a^*$. The sign of those terms depends on the sign of $z_x\cdots z_1$; so there are $k^2l(2l-1)+kl(k-1)(2k-2)$ positive terms and $k^2l(2k-2)+kl(k-1)(2l-1)$ negative terms in total. So if the derivation is locally nilpotent of degree three for $a^*$, then these kinds of terms need to cancel each other out, so $2(k-l)=1$, which is impossible.

\item 
Suppose now $\theta_n^{t}(a^*)$ is not $0$ and its terms of minimal degree in $a$ and $a^*$ contain $tk-t$ times the arrow $a$ and $tl-(t-1)$ times the arrow $a^*$. Suppose furthermore that $\theta_n^{t}(a^*)$ contains $s_+$ positive terms of minimal degree and $s_-$ negative ones and $s_+\neq s_-$. 

\item Now, terms of the lowest degree of $\theta_n^{t+1}(a^*)$ are those of the form 
$$z_x\cdots \theta_n(a^*)\cdots z_1\qquadtext{ and }z_x\cdots\theta_n(a)\cdots z_1$$
with $z_x\cdots z_1$ being a term of minimal degree in $a$ and $a^*$ of $\theta_n^{t}(a^*)$. Both above terms contain $(t+1)k-(t+1)$ times the arrow $a$ and $(t+1)l-t$ times the arrow $a^*$. The sign of those terms depends on the sign of $z_x\cdots z_1$; so there are $s_+(tl-(t-1))+s_-(tk-t)$ positive terms and $s_+(tk-t)+s_-(tl-(t-1))$ negative terms in total. So if the derivation is locally nilpotent of degree $t+1$ for $a^*$, then these kinds of terms need to cancel each other out, so $t(k-l)=1$, which is impossible.
\end{itemize}

Thus, the derivation $\theta_n$ is not locally nilpotent, because for any $t\in\N$, $\theta_n^t(a^*)\neq 0$. 
\end{proof}

A special case of the previous lemma is the following:

\begin{lemma}
If a necklace $n$ contains only $*$-ed or only un-$*$-ed arrows (i.e.\ $\Supp(n)\subseteq L$ or $\Supp(n)\subseteq L^*$) then the derivation $\theta_n$ is locally nilpotent.
\end{lemma}

\begin{definition}
Using a locally nilpotent derivation $\theta_n$, we can define a \emph{symplectic flow} $\gamma_n$ on $\C{\Delta}$ by $$\gamma_{n, \rho} \defeq e^{\rho\theta_n}=\sum_{j=0}^\infty \frac{\rho^j}{j!}\theta_n^j.$$
\end{definition}

Given a double quiver $\Delta$, we're interested in the group generated by all symplectic flows on $\Delta$ coming from locally nilpotent derivations,
$$\Flow^\Delta\defeq\gend{\gamma_{n,\rho}}{n \text{ a one-way necklace}, \rho\in\C}.$$
Let $\set{n_1,\ldots,n_l}$ be all primitive one-way necklaces. We would like for $m_1,m_2\in\C[n_1,\ldots,n_l]$ and $\rho_1,\rho_2\in\C$ that $$\gamma_{m_2,\rho_2}\circ\gamma_{m_1,\rho_1}=\gamma_{m,\rho}$$ for a certain $m\in\C[n_1,\ldots,n_l]$ and $\rho\in\C$. Unfortunately, this is not the case in general.  This can be easily seen in the case of the Calogero-Moser quiver from example~\ref{CMquiver}, with necklaces $m_1=\frac{1}{k+1}b^{k+1}$ and $m_2=\frac{1}{l+1}{b^*}^{l+1}$. The statement does hold for the certain (abelian) subgroups of the form
$$\Flow^\Delta_{n_1,\ldots,n_k}\defeq\gend{\gamma_{m,\rho}}{m\in\C[n_1,\ldots,n_k], \rho\in\C}\subseteq \Flow^\Delta,$$
with $k\leq l$.

\begin{lemma}\label{flowgen}
Let $n_1,\ldots,n_k$ be a set of primitive one-way necklaces in a double quiver $\Delta$ and for which
$$\forall a\in\arr: \qquad \set{a,a^*}\not\subseteq \bigcup_{i=1}^k\Supp(n_i)$$
then
$$\Flow^\Delta_{n_1,\ldots,n_k}=\setd{\gamma_{m,\rho}}{m\in\C[n_1,\ldots,n_k], \rho\in\C}$$
with $\gamma_{m_2,\rho_2}\circ\gamma_{m_1,\rho_1}=\gamma_{m,1}$ where $m=\rho_2 m_2+\rho_1 m_1$.
\end{lemma}

\begin{proof}
Proof left to the reader.
\end{proof}

 We will briefly put forward three examples: one related to Calogero-Moser flows, one related to the projective modular group and one related to $\GL_2(\Z)$.

\begin{example}\label{CMquiver}
Consider the symplectic Calogero-Moser quiver,
$$
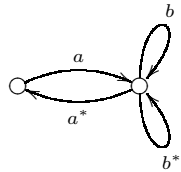
\xymatrix@=4em{\vtx{} \ar@/^/[r]^a & 
\vtx{} \ar@/^/[l]^{a^*} \ar@(u,ru)[]^{b} \ar@(d,rd)[]_{b^*}}
$$
We can define a symplectic structure on this quiver by taking $L=\set{a,b}$. Under composition, all necklaces are now generated by the elementary necklaces $aa^*$, $b$ and $b^*$. All Lie brackets between those elementary necklaces are zero.

The only necklaces in $\Delta$ containing only un-$*$-ed arrows are monomials in $b$. The symplectic derivation $\theta_n$ corresponding to the necklace $n=\frac{1}{k+1}b^{k+1}$ is defined by
$$\theta_n(a)=0 \qquad \theta_n(a^*)=0 \qquad \theta_n(b)=0 \qquad \theta_n(b^*)=b^k,$$
and therefore,
$$\theta_n^{m+1}((b^*)^m)=0,$$
so the corresponding flow $\gamma_{n,\rho}$ is well defined, by
$$\gamma_{n,\rho}(a)=a \qquad \gamma_{n,\rho}(a^*)=a^* \qquad \gamma_{n,\rho}(b)=b \qquad \gamma_{n,\rho}(b^*)=b^*+\rho b^k.$$
In this case, the moment element $m$ turns out to be
$$m=[a,a^*]+[b,b^*]=-a^*ae_1+(aa^*+bb^*-b^*b)e_2.$$
We can now take as dimension vector $\alpha=(1,k)$ and look at the \emph{deformed preprojective algebra} (\cite{LBBocklandt}), 
$$\Pi_\Lambda=\frac{\C{\Delta}}{(m-\lambda)}$$
where $\lambda=-ke_1+1e_2$. Then
$$\Iss_\alpha\Pi_\lambda=\Rep_\alpha\Pi_\lambda \aquot \GL_\alpha = \Calo_k,$$
the \emph{Calogero-Moser phase space}. This space $\Calo_k$ can be identified with pairs of $k\times k$-matrices
$$(X,Z)\in\M_k(\C)\times \M_k(\C) \quadtext{ where } \rank\left([X,Z]+I_k\right)=1.$$
For more details on the Calogero-Moser phase space, see \cite{Wil98} and the link with necklace Lie algebras is investigated in \cite{LBBocklandt}. These matrices $X$ and $Z$ correspond to the loops $b$ and $b^*$ in $\Delta$. As the coordinate ring of the quotient variety $\Iss_\alpha \Delta$ (and thus the quotient $\Iss_\alpha \Pi_\lambda$) is generated by traces along oriented  cycles in $\Delta$, we see that the automorphism $\gamma_{n,\rho}$ defines a $\GL_n$-invariant flow on $\Rep_\alpha\Pi_\lambda$ defined by
$$(X,Z,u,v)\mapsto(X+tZ^n,Z,u,v)$$
which generate the \emph{Calogero-Moser flows} on $\Calo_k$ as defined in \cite{Wil98}.

For the Calogero-Moser quiver, there are two abelian subgroups $\Flow^\Delta_n$ arising from lemma~\ref{flowgen}: one for $n=b$ and one for $n=b^*$.
\end{example}

\begin{example}\label{ex:liepsl2}
We can define a symplectic structure on the quiver $\Delta$ corresponding to $\C\PSL_2(\Z)$, 
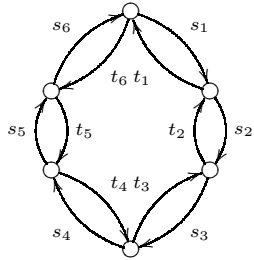
$$\xymatrix{
& \vtx{} \ar@/^/[dl]^{t_6} \ar@/^/[dr]^{s_1} & \\
\vtx{}  \ar@/^/[ur]^{s_6} \ar@/^/[d]^{t_5} & & \vtx{}  \ar@/^/[ul]^{t_1} \ar@/^/[d]^{s_2} \\
\vtx{}  \ar@/^/[dr]^{t_4} \ar@/^/[u]^{s_5} & & \vtx{}  \ar@/^/[dl]^{s_3} \ar@/^/[u]^{t_2} \\
& \vtx{} \ar@/^/[ur]^{t_3} \ar@/^/[ul]^{s_4} 
}$$
by taking $L=\set{s_1,\ldots,s_6}$ and $s_i^*=t_i$. Under compositition, all necklaces are now generated by the elementary necklaces:
$$L_i=s_it_i \qquad S=s_6s_5s_4s_3s_2s_1 \qquad T=t_1t_2t_3t_4t_5t_6$$
All Lie brackets between those elementary necklaces are
$$[L_i,L_j]=0 \qquad [L_i,S]=S \qquad [T,L_i]=T$$
$$[S,T]=
\sum_{i=1}^6 \left(\prod_{j\neq i}L_j\right)
=\sum_{i=1}^6\left( 
\centerprent{\xymatrix@=1em{
& \vtx{} \ar@/^.2em/[dl] \ar@{}[dr]|{L_i} & \\
\vtx{}  \ar@/^.2em/[ur] \ar@/^.2em/[d] & & \vtx{}  \ar@/^.2em/[d] \\
\vtx{}  \ar@/^.2em/[dr] \ar@/^.2em/[u] & & \vtx{}  \ar@/^.2em/[dl] \ar@/^.2em/[u] \\
& \vtx{} \ar@/^.2em/[ur] \ar@/^.2em/[ul]
}} 
\right)
$$

The only one-way necklaces in $\Delta$ are monomials in $S$ and monomials in $T$. The symplectic derivation $\theta_n$ corresponding to the necklace $n=\frac{1}{k+1}S^{k+1}$ is defined by
$$\theta_n(s_i)=0 \qquad \theta_n(t_i)=  s_{i-1} \cdots s_1 S^k s_6 \cdots s_{i+1},$$
and therefore,
$\theta_n^{6d+1}(T^d)=0$,
so the corresponding flow $\gamma_{n,\rho}$ is well defined, by
$$\gamma_{n,\rho}(s_i)=s_i \qquad \gamma_{n,\rho}(t_i)=t_i+\rho(s_{i-1} \cdots s_1 S^k s_6 \cdots s_{i+1}).$$
For $\PSL_2(\Z)$, there are two abelian subgroups $\Flow^\Delta_n$ arising from lemma~\ref{flowgen}: one for $n=S$ and one for $n=T$.
\end{example}

\begin{example}\label{ex:liegl2}
Let us name the arrows in the larger component of the local quiver of $\GL_2(\Z)$ (the smaller component is not double),
$$
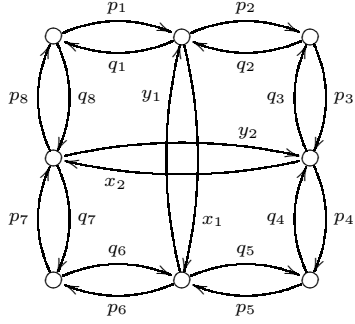
\centerprent{
\xymatrix@=4em{
\vtx{~} \ar@/^/[r]^{p_1} \ar@/^/[d]^{q_8} &
\vtx{~} \ar@/^/[l]^{q_1} \ar@/^/[r]^{p_2} \ar@/^/[dd]^(.75){x_1} &
\vtx{~} \ar@/^/[l]^{q_2} \ar@/^/[d]^{p_3} \\
\vtx{~} \ar@/^/[u]^{p_8} \ar@/^/[rr]^(.75){y_2} \ar@/^/[d]^{q_7} & &
\vtx{~} \ar@/^/[ll]^(.75){x_2} \ar@/^/[u]^{q_3} \ar@/^/[d]^{p_4} \\
\vtx{~} \ar@/^/[u]^{p_7} \ar@/^/[r]^{q_6} &
\vtx{~} \ar@/^/[uu]^(.75){y_1} \ar@/^/[l]^{p_6} \ar@/^/[r]^{q_5} &
\vtx{~} \ar@/^/[l]^{p_5} \ar@/^/[u]^{q_4}
}}$$

Let $L=\set{p_1,\ldots,p_8,x_1,x_2}$, and let $p_i^*=q_i$ and $x_i^*=y_i$. 
Under composition, all necklaces are now generated by the elementary necklaces: 
$$
\centerprent{\xymatrix@=1em{
\vtx{} \ar@/^.2em/[r] &
\vtx{} \ar@/^.2em/[r] &
\vtx{} \ar@/^.2em/[d] \\
\vtx{} \ar@/^.2em/[u] & &
\vtx{} \ar@/^.2em/[d] \\
\vtx{} \ar@/^.2em/[u] &
\vtx{} \ar@/^.2em/[l] &
\vtx{} \ar@/^.2em/[l]
}}
=G^p=p_8\cdots p_1
\qquad \qquad 
\centerprent{\xymatrix@=1em{
\vtx{} \ar@/^.2em/[d] &
\vtx{} \ar@/^.2em/[l] &
\vtx{} \ar@/^.2em/[l] \\
\vtx{} \ar@/^.2em/[d] & &
\vtx{} \ar@/^.2em/[u] \\
\vtx{} \ar@/^.2em/[r] &
\vtx{} \ar@/^.2em/[r] &
\vtx{} \ar@/^.2em/[u]
}} 
=G^q=p_1\cdots p_8
$$

$$
\centerprent{\xymatrix@=1em{
\vtx{} \ar@/^.2em/[r] & 
\vtx{} \ar@/^.2em/[l] & {}\\
\\
{}
}}
=L_i=p_iq_i
\qquad \qquad
\centerprent{\xymatrix@=1em{
\\
\vtx{} \ar@/^.2em/[rr] & &
\vtx{} \ar@/^.2em/[ll] \\ {}
}}
=M_i=x_iy_i
$$

$$
\centerprent{\xymatrix@=1em{
\vtx{} \ar@/^.2em/[r] &
\vtx{} \ar@/^.2em/[r] &
\vtx{} \ar@/^.2em/[d] \\
\vtx{} \ar@/^.2em/[u] & &
\vtx{} \ar@/^.2em/[ll] \\
{}
}} 
=\left\{
\begin{tabular}{C C}
K^p_{x_1}=x_1p_1p_8p_7p_6 
&
K^p_{x_2}=x_2p_3p_2p_1p_8 
\\
K^p_{y_1}=y_1p_5p_4p_3p_2
&
K^p_{y_2}=y_2p_7p_6p_5p_4
\end{tabular}
\right. 
$$

$$
\centerprent{\xymatrix@=1em{
\vtx{} \ar@/^.2em/[d] &
\vtx{} \ar@/^.2em/[l] &
\vtx{} \ar@/^.2em/[l] \\
\vtx{} \ar@/^.2em/[rr] & &
\vtx{} \ar@/^.2em/[u] \\
{}
}} 
=\left\{
\begin{tabular}{C C}
K^q_{x_1}=x_1q_2q_3q_4q_5
&
K^q_{x_2}=x_2q_4q_5q_6q_7
\\
K^q_{y_1}=y_1q_6q_7q_8q_1
&
K^q_{y_2}=y_2q_8q_1q_2q_3
\end{tabular}
\right. 
$$

$$
\centerprent{\xymatrix@=1em{
\vtx{} \ar@/^.2em/[r] &
\vtx{} \ar@/^.2em/[dd] &
\\
\vtx{} \ar@/^.2em/[u] & &
\vtx{} \ar@/^.2em/[ll] \\
&
\vtx{} \ar@/^.2em/[r] &
\vtx{} \ar@/^.2em/[u]
}} 
=\left\{
\begin{tabular}{C C}
E_{x_1x_2}=p_1p_8x_2q_4q_5x_1
&
E_{x_1y_2}=q_2q_3y_2p_7p_6x_1
\\
E_{y_1x_2}=q_6q_7x_2p_3p_2y_1
&
E_{y_1y_2}=p_5p_4y_2q_8q_1y_1
\end{tabular}
\right. 
$$

We calculate the brackets between those elementary necklaces, and where the occasional $\centerprent{\xymatrix@=1em{
\vtx{} & &
\vtx{} \ar@/^.2em/[ll]
}}$ in the first part of the bracket is an $x_i$. Some interesting non-zero brackets are: 
$$
\left[
\centerprent{\xymatrix@=1em{
\vtx{} \ar@/^.2em/[r] &
\vtx{} \ar@/^.2em/[r] &
\vtx{} \ar@/^.2em/[d] \\
\vtx{} \ar@/^.2em/[u] & &
\vtx{} \ar@/^.2em/[ll] \\
{}
}} 
\ ,\ 
\centerprent{\xymatrix@=1em{
\\
\vtx{} \ar@/^.2em/[rr] & &
\vtx{} \ar@/^.2em/[d] \\
\vtx{} \ar@/^.2em/[u] &
\vtx{} \ar@/^.2em/[l] &
\vtx{} \ar@/^.2em/[l] 
{}
}} 
\right]
=
\centerprent{\xymatrix@=1em{
\vtx{} \ar@/^.2em/[r] &
\vtx{} \ar@/^.2em/[r] &
\vtx{} \ar@/^.2em/[d] \\
\vtx{} \ar@/^.2em/[u] & &
\vtx{} \ar@/^.2em/[d] \\
\vtx{} \ar@/^.2em/[u] &
\vtx{} \ar@/^.2em/[l] &
\vtx{} \ar@/^.2em/[l]
}}
$$


$$
\left[
\centerprent{\xymatrix@=1em{
\vtx{} \ar@/^.2em/[r] &
\vtx{} \ar@/^.2em/[r] &
\vtx{} \ar@/^.2em/[d] \\
\vtx{} \ar@/^.2em/[u] & &
\vtx{} \ar@/^.2em/[d] \\
\vtx{} \ar@/^.2em/[u] &
\vtx{} \ar@/^.2em/[l] &
\vtx{} \ar@/^.2em/[l]
}}
\ ,\ 
\centerprent{\xymatrix@=1em{
\vtx{} \ar@/^.2em/[d] &
\vtx{} \ar@/^.2em/[l] &
\vtx{} \ar@/^.2em/[l] \\
\vtx{} \ar@/^.2em/[rr] & &
\vtx{} \ar@/^.2em/[u] \\
{}
}} 
\right]
=
\centerprent{\xymatrix@=1em{
\vtx{} \ar@/^.2em/[r] &
\vtx{} \ar@/^.2em/[l] \ar@/^.2em/[r] &
\vtx{} \ar@/^.2em/[l] \ar@/^.2em/[d] \\
\vtx{} \ar@/^.2em/[rr] & &
\vtx{} \ar@/^.2em/[u] \ar@/^.2em/[d] \\
\vtx{} \ar@/^.2em/[u] &
\vtx{} \ar@/^.2em/[l] &
\vtx{} \ar@/^.2em/[l] 
}}
+
\centerprent{\xymatrix@=1em{
\vtx{} \ar@/^.2em/[d] &
\vtx{} \ar@/^.2em/[r] &
\vtx{} \ar@/^.2em/[l] \ar@/^.2em/[d] \\
\vtx{} \ar@/^.2em/[u] \ar@/^.2em/[rr] & &
\vtx{} \ar@/^.2em/[u] \ar@/^.2em/[d] \\
\vtx{} \ar@/^.2em/[u] &
\vtx{} \ar@/^.2em/[l] &
\vtx{} \ar@/^.2em/[l] 
}}
+
\centerprent{\xymatrix@=1em{
\vtx{} \ar@/^.2em/[r] \ar@/^.2em/[d] &
\vtx{} \ar@/^.2em/[l] &
\vtx{} \ar@/^.2em/[d] \\
\vtx{} \ar@/^.2em/[u] \ar@/^.2em/[rr] & &
\vtx{} \ar@/^.2em/[u] \ar@/^.2em/[d] \\
\vtx{} \ar@/^.2em/[u] &
\vtx{} \ar@/^.2em/[l] &
\vtx{} \ar@/^.2em/[l] 
}}
+
\centerprent{\xymatrix@=1em{
\vtx{} \ar@/^.2em/[r] \ar@/^.2em/[d] &
\vtx{} \ar@/^.2em/[l] \ar@/^.2em/[r] &
\vtx{} \ar@/^.2em/[l] \\
\vtx{} \ar@/^.2em/[u] \ar@/^.2em/[rr] & &
\vtx{} \ar@/^.2em/[d] \\
\vtx{} \ar@/^.2em/[u] &
\vtx{} \ar@/^.2em/[l] &
\vtx{} \ar@/^.2em/[l] 
}}
$$
$$
\left[
\centerprent{\xymatrix@=1em{
\vtx{} \ar@/^.2em/[r] &
\vtx{} \ar@/^.2em/[r] &
\vtx{} \ar@/^.2em/[d] \\
\vtx{} \ar@/^.2em/[u] & &
\vtx{} \ar@/^.2em/[ll] \\
{}
}} 
\ ,\ 
\centerprent{\xymatrix@=1em{
\vtx{} \ar@/^.2em/[d] &
\vtx{} \ar@/^.2em/[l] &
\\
\vtx{} \ar@/^.2em/[d] & &
\\
\vtx{} \ar@/^.2em/[r] &
\vtx{} \ar@/^.2em/[uu] &
}} 
\right]
=
\centerprent{\xymatrix@=1em{
\vtx{}  \ar@/^.2em/[d] &
\vtx{}  \ar@/^.2em/[r] &
\vtx{} \ar@/^.2em/[d] \\
\vtx{} \ar@/^.2em/[u]  \ar@/^.2em/[d] & &
\vtx{} \ar@/^.2em/[ll]  \\
\vtx{}  \ar@/^.2em/[r] &
\vtx{} \ar@/^.2em/[uu] 
}}
+
\centerprent{\xymatrix@=1em{
\vtx{} \ar@/^.2em/[r]  &
\vtx{} \ar@/^.2em/[l] \ar@/^.2em/[r] &
\vtx{} \ar@/^.2em/[d] \\
\vtx{}  \ar@/^.2em/[d] & &
\vtx{} \ar@/^.2em/[ll]  \\
\vtx{}  \ar@/^.2em/[r] &
\vtx{} \ar@/^.2em/[uu] 
}}
$$

$$
\left[
\centerprent{\xymatrix@=1em{
\vtx{} \ar@/^.2em/[r] &
\vtx{} \ar@/^.2em/[dd] &
\\
\vtx{} \ar@/^.2em/[u] & &
\vtx{} \ar@/^.2em/[ll] \\
&
\vtx{} \ar@/^.2em/[r] &
\vtx{} \ar@/^.2em/[u]
}} 
\ ,\ 
\centerprent{\xymatrix@=1em{
&
\vtx{} \ar@/^.2em/[r] &
\vtx{} \ar@/^.2em/[d]
\\
\vtx{} \ar@/^.2em/[d] & &
\vtx{} \ar@/^.2em/[ll] \\
\vtx{} \ar@/^.2em/[r] &
\vtx{} \ar@/^.2em/[uu] &
}} 
\right]
=
\pm
\centerprent{\xymatrix@=1em{
\vtx{} \ar@/^.2em/[r] &
\vtx{} \ar@/^.2em/[r] &
\vtx{} \ar@/^.2em/[d] \\
\vtx{} \ar@/^.2em/[u] \ar@/^.2em/[d] & &
\vtx{} \ar@/^.2em/[ll] \ar@/^.4em/[ll] \\
\vtx{} \ar@/^.2em/[r] & 
\vtx{} \ar@/^.2em/[r] &
\vtx{} \ar@/^.2em/[u]
}}
$$

The necklaces in $\Delta$ containing only $*$-ed arrows are generated by 
$$\set{G^p, K^p_{x_1}, K^p_{x_2}}.$$
Let's consider the three derivations corresponding to these elementary necklaces. First, for $n=G^p$, we have
$$\theta_n(p_i)=\theta_n(x_i)=\theta_n(y_i)=0 \qquad \theta_n(q_i)=p_{i-1}\cdots p_1p_8\cdots p_{i+1}$$
giving rise to the flow $\gamma_{n,\rho}$,
$$\gamma_{n,\rho}(p_i)=p_i \qquad \gamma_{n,\rho}(x_i)=x_i \qquad \gamma_{n,\rho}(y_i)=y_i$$
$$\gamma_{n,\rho}(q_i)=q_i+\rho(p_{i-1}\cdots p_1p_8\cdots p_{i+1})$$
For the necklace $n=K^p_{x_1}$, the derivation becomes,
$$\theta_n(p_i)=\theta_n(x_2)=\theta_n(y_i)=0$$
$$\theta_n(q_i)=p_{i-1}\cdots p_6x_1p_1\cdots p_{i+1} \quad\text{for }i=1,6,7,8$$
$$\theta_n(q_i)=0 \quad\text{for }i=2,3,4,5$$
giving rise to the flow $\gamma_{n,\rho}$,
$$\gamma_{n,\rho}(p_i)=p_i \qquad \gamma_{n,\rho}(x_2)=x_2 \qquad \gamma_{n,\rho}(y_i)=y_i$$
$$\gamma_{n,\rho}(q_i)=q_i+\rho(p_{i-1}\cdots p_6x_1p_1\cdots p_{i+1}) \quad\text{for }i=1,6,7,8$$
$$\gamma_{n,\rho}(q_i)=q_i \quad\text{for }i=2,3,4,5$$
Finally and analogously, we find for the necklace $n=K^p_{x_2}$,
$$\theta_n(p_i)=\theta_n(x_1)=\theta_n(y_i)=0$$
$$\theta_n(q_i)=p_{i-1}\cdots p_8x_2p_3\cdots p_{i+1} \quad\text{for }i=1,2,3,8$$
$$\theta_n(q_i)=0 \quad\text{for }i=4,5,6,7$$
giving rise to the flow $\gamma_{n,\rho}$,
$$\gamma_{n,\rho}(p_i)=p_i \qquad \gamma_{n,\rho}(x_1)=x_1 \qquad \gamma_{n,\rho}(y_i)=y_i$$
$$\gamma_{n,\rho}(q_i)=q_i+\rho(p_{i-1}\cdots p_8x_2p_3\cdots p_{i+1}) \quad\text{for }i=1,2,3,8$$
$$\gamma_{n,\rho}(q_i)=q_i \quad\text{for }i=4,5,6,7$$

For $\GL_2(\Z)$, there are sixteen abelian subgroups $\Flow^\Delta_{n_1,\ldots,n_k}$ arising from lemma~\ref{flowgen}. There are eight groups corresponding to
$$
n_1=
\centerprent
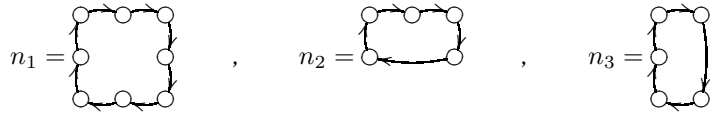
{\xymatrix@=1em{
\vtx{} \ar@/^.2em/[r] &
\vtx{} \ar@/^.2em/[r] &
\vtx{} \ar@/^.2em/[d] \\
\vtx{} \ar@/^.2em/[u] & &
\vtx{} \ar@/^.2em/[d] \\
\vtx{} \ar@/^.2em/[u] &
\vtx{} \ar@/^.2em/[l] &
\vtx{} \ar@/^.2em/[l]
}}
\qquadtext{ , }
n_2=
\centerprent{\xymatrix@=1em{
\vtx{} \ar@/^.2em/[r] &
\vtx{} \ar@/^.2em/[r] &
\vtx{} \ar@/^.2em/[d] \\
\vtx{} \ar@/^.2em/[u] & &
\vtx{} \ar@/^.2em/[ll] \\
&
&
}}
\qquadtext{ , }
n_3=
\centerprent{\xymatrix@=1em{
\vtx{} \ar@/^.2em/[r] &
\vtx{} \ar@/^.2em/[dd] &
\\
\vtx{} \ar@/^.2em/[u] &  &
\\
\vtx{} \ar@/^.2em/[u] &
\vtx{} \ar@/^.2em/[l] &
}} 
$$
and eight corresponding to
$$
n_1=
\centerprent
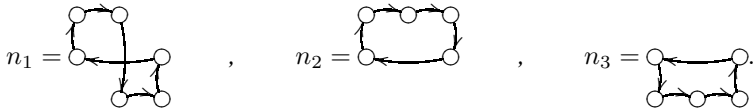
{\xymatrix@=1em{
\vtx{} \ar@/^.2em/[r] &
\vtx{} \ar@/^.2em/[dd] &
\\
\vtx{} \ar@/^.2em/[u] & &
\vtx{} \ar@/^.2em/[ll] \\
&
\vtx{} \ar@/^.2em/[r] &
\vtx{} \ar@/^.2em/[u]
}}
\qquadtext{ , }
n_2=
\centerprent{\xymatrix@=1em{
\vtx{} \ar@/^.2em/[r] &
\vtx{} \ar@/^.2em/[r] &
\vtx{} \ar@/^.2em/[d] \\
\vtx{} \ar@/^.2em/[u] & &
\vtx{} \ar@/^.2em/[ll] \\
&
&
}}
\qquadtext{ , }
n_3=
\centerprent{\xymatrix@=1em{
&
&
\\
\vtx{} \ar@/^.2em/[d] & &
\vtx{} \ar@/^.2em/[ll] \\
\vtx{} \ar@/^.2em/[r] &
\vtx{} \ar@/^.2em/[r] &
\vtx{} \ar@/^.2em/[u] 
}}
.$$

\end{example}

\section{Compactifications}
 
 In the foregoing section we have seen that one-way necklaces induce flows on the affine quotient varieties $\wis{iss}_{\alpha}~\ST$. In order to turn these flows into integral curves it would be interesting to have natural compactifications of these quotient varieties. In this section we will show that this is indeed possible following the suggestion from \cite{LBncm} to construct non-commutative compact manifolds.
 
 Let $\ST$ be a tree of semi-simple algebras with associated Zariski quiver $\Gamma = \Gamma(\ST)$. In section~\ref{Zariski} we have seen that we can identify isomorphism classes of finite dimensional $\ST$-representations with those of a certain universal localization $\C \Gamma_{\Sigma}$ of the path algebra of this quiver without oriented cycles.
 
\begin{definition} A dimension vector $\alpha = (\alpha_i^{(v)} : v \in V, 1 \leq i \leq n_v)$ of $\Gamma$ is said to be a {\em weak $n$-dimension vector} provided for all $v \in V$ we have
\[
\sum_{i=1}^{n_v} d_v(i) \alpha_i^{(v)} = n \]
 For any edge $\xymatrix{\vtx{v} \ar@{-}[r] & \vtx{w}}$ in $T$ we denote with $\Gamma_e$ the bipartite subquiver of $\Gamma = \Gamma(\ST)$ defined by the embeddings $S_v \lInto S_e \rInto S_w$.
 
 Given a weak $n$-dimension vector $\alpha$ and an edge $e \in T$ and $M \in \wis{rep}_{\alpha}~\Gamma$ we say that
\begin{itemize}
\item{$M$ is {\em $e$-semistable} if and only if for every subrepresentation $N$ of $M | \Gamma_e$ of dimension vector $(k_1,\hdots,k_{n_v},l_1,\hdots,l_{n_w})$ we have
\[
\sum_{i=1}^{n_w} l_i d_w(i) \geq \sum_{i=1}^{n_v} k_i d_v(i) \]}
  \item{$M$ is {\em $e$-stable} if and only if for every subrepresentation $N$ of $M | \Gamma_e$ of dimension vector $(k_1,\hdots,k_{n_v},l_1,\hdots,l_{n_w})$ we have
\[
\sum_{i=1}^{n_w} l_i d_w(i) > \sum_{i=1}^{n_v} k_i d_v(i) \]}
\item{$M$ is $\ST$-semistable (resp. $\ST$-stable) if and only if $M$ is $e$-semistable (resp. $e$-stable) for every edge $e \in T$.}
\end{itemize}
\end{definition}
 
 We recall from \cite[Prop 8]{LBqaq} that every $\Gamma$-representation corresponding to a finite dimensional $\ST$-representation is $e$-semistable and is $e$-stable if and only if the $\ST$-representation is simple. That is, $\wis{rep}~\ST$ is an Abelian subcategory of the Abelian category $\wis{sst}~\Gamma$ of all $\ST$-semistable representations of $\Gamma$ (for varying weak $n$-dimension vectors). Observe that by additivity of dimension vectors in short exact sequences, kernels and cokernels of morphisms between two $\ST$-stable representations are again $\ST$-semistable. Modifying the argument of \cite{LBqaq} one can cover $\wis{sst}~\Gamma$ with subcategories isomorphic to $\wis{rep}~\C \Gamma_{S}$, the finite dimensional representations of suitable universal localizations of the finite dimensional path algebra $\C \Gamma$. We outline the main ideas.
 
 Fix an edge $\xymatrix{\vtx{v} \ar@{-}[r] & \vtx{w}}$ in $T$, let $\Gamma_e$ be the corresponding bipartite quiver and $\alpha_e = (a_1,\hdots,a_{n_v},b_1,\hdots,b_{n_w})$ the restriction of $\alpha$ to $\Gamma_e$. Consider the character
$ \theta_e~:~GL(\alpha_e) \rTo \C^*$ sending $(g_1,\hdots,g_{n_v},h_1,\hdots,h_{n_w})$ to
\[
 det(g_1)^{d_v(1)} \hdots det(g_{n_v})^{d_v(n_v)} det(h_1)^{-d_w(1)} \hdots det(h_{n_w})^{-d_w(n_w)} \]
 A polynomial function $f \in \C[\wis{rep}_{\alpha_e}~\Gamma_e]$ is said to be $\theta_e$-semi-invariant of weight $l \in \N$ whenever
\[
 g.f = \theta_e(g)^l f \qquad \forall g \in GL(\alpha_e) \]
 The weight function turns the ring of polynomial $\theta_e$-semi-invariants into a positively graded algebra $\C[\wis{rep}_{\alpha_e}~\Gamma_e]^{GL(\alpha_e),\theta_e}$ with part of degree zero reduced to $\C$ (as $\Gamma_e$ is bipartite and in particular does not have oriented cycles). Hence, it defines a projective variety
\[
\wis{moduli}^{\theta_e}_{\alpha_e}~\Gamma_e = \wis{proj}~\C[\wis{rep}_{\alpha_e}~\Gamma_e]^{GL(\alpha_e),\theta_e} \]
 called the {\em moduli space of $e$-semistable $\alpha_e$-dimensional representations of $\Gamma_e$}. From \cite{King} we recall that the $\C$-points of $\wis{moduli}^{\theta_e}_{\alpha_e}~\Gamma_e$ correspond to isomorphism classes of direct sums of $e$-stable representations of $\Gamma_e$ of total dimension $\alpha_e$. 
 
 Repeating this construction for all edges $\{ e_1,\hdots,e_z \}$ we therefore obtain an immersion of the quotient variety $\wis{iss}_{\alpha}~\ST$ into a product of projective varieties
\[
\wis{iss}_{\alpha}~\ST \rInto \wis{moduli}^{\theta_{e_1}}_{\alpha_{e_1}}~\Gamma_{e_1} \times \hdots \times \wis{moduli}^{\theta_{e_z}}_{\alpha_{e_z}}~\Gamma_{e_z} \]
 In order to study these varieties as well as to describe the link with universal localizations of $\C \Gamma$ we recall the generating result of $\theta_e$-semi-invariants by {\em determinantal semi-invariants} due to Aidan Schofield and Michel Van den Bergh, \cite{SchofVdB}.
 
 For a fixed edge $\xymatrix{\vtx{v} \ar@{-}[r] & \vtx{w}}$ in $T$ and weight $l_e \in \N_+$ consider a block-matrix $\Delta_e$ with entries in $\C \Gamma_e$
\[
\Delta_e =  \begin{bmatrix} 
\begin{array}{c | c c | c}
A_{11} & \hdots & \hdots & A_{n_v1} \\
& & & \\
\hline & & & \\
\vdots & & & \vdots \\
& & & \\
\hline 
& & & \\
A_{1n_w} & \hdots & \hdots & A_{n_v n_w} \\
& & & 
\end{array}
\end{bmatrix}
\]
where the block-matrix $A_{ij}$ has size $l_e d_w(j) \times l_e d_v(i)$ and all its entries are linear combinations of arrows in $\Gamma_e$ from vertex $\mu_i^{(v)}$ to vertex $\mu_j^{(w)}$.

Evaluating the matrix $\Delta_e$ in $V \in \wis{rep}_{\alpha_e}~\Gamma_e$ for a weak $n$-dimension vector $\alpha$ we obtain a square matrix and its determinant
\[
det~\Delta_e(-)~:~\wis{rep}_{\alpha_e}~\Gamma_e \rTo \C \]
is a $\theta_e$-semi-invariant polynomial function of weight $l$. In \cite{SchofVdB} it is proved that these {\em determinantal semi-invariants} generate the ring of semi-invariants $\C[\wis{rep}_{\alpha_e}~\Gamma_e]^{GL(\alpha_e),\theta_e}$. 

Observe, as before, that the matrix $\Delta_e$ defines a morphism of right $\C \Gamma$-modules between finitely generated projective $\C \Gamma$-modules.
For $l = (l_{e_1},\hdots,l_{e_z}) \in \N^z_+$ and $\Delta = \{ \Delta_{e_1},\hdots,\Delta_{e_z} \}$ consider the universal localization $\C \Gamma_{\Delta}$. Observe that the algebra $\C \Gamma_{\Sigma}$ constructed in section~\ref{Zariski} having the same finite dimensional representations as $\ST$ is a special case of this construction. Similarly, for each $\Delta$, the corresponding quotient variety
\[
\mathbb{X}_{\alpha}(\Delta) = \wis{iss}_{\alpha}~\C \Gamma_{\Delta} \]
is an affine open subset of the projective variety $ \wis{moduli}^{\theta_{e_1}}_{\alpha_{e_1}}~\Gamma_{e_1} \times \hdots \times \wis{moduli}^{\theta_{e_z}}_{\alpha_{e_z}}~\Gamma_{e_z}$. By the generating result, these affine open subsets $\mathbb{X}_{\alpha}(\Delta)$ form a Zariski cover. Phrased in categorical terms we have

\begin{theorem} The abelian category $\wis{sst}~\Gamma$ of all $\ST$-semistable finite dimensional representations of the Zariski quiver $\Gamma = \Gamma(\ST)$ is the union
\[
\wis{sst}~\Gamma = \bigcup_{\Delta}~\wis{rep}~\C \Gamma_{\Delta} \]
where $\Delta$ runs over all $z$-tuples of matrices constructed above. Observe that all $\C \Gamma_{\Delta}$ are qurves so $\wis{ss}~\Gamma$ can be viewed as a non-commutative compact manifold.
\end{theorem}

To perform a local study of the projective varieties which compactify $\wis{iss}_{\alpha}~\ST$ we need to be able to determine the \'etale quiver of the universal localizations $\C \Gamma_{\Delta}$. Again, the vertices correspond to the semigroup generators of a certain sub-semigroup $\mathcal{S}_{\Delta}$ of all dimension vectors of $\Gamma$ and the Euler-form $\chi_{\Gamma}$ can then be used to determine the arrows in the \'etale quiver. 

\begin{theorem} For fixed $l=(l_{e_1},\hdots,l_{e_z}) \in \N^z_+$ and $\Delta = \{ \Delta_{e_1},\hdots,\Delta_{e_z} \}$, the semigroup $\mathcal{S}_{\Delta}$ consists of the weak $n$-dimension vectors $\alpha$ such that for every edge $\xymatrix{\vtx{v} \ar@{-}[r] & \vtx{w}}$ in $T$, weight $l_e \in \N_+$ and block-matrix $\Delta_e$
\[
\Delta_e =  \begin{bmatrix} 
\begin{array}{c | c c | c}
A_{11} & \hdots & \hdots & A_{n_v1} \\
& & & \\
\hline & & & \\
\vdots & & & \vdots \\
& & & \\
\hline & & & \\
A_{1n_w} & \hdots & \hdots & A_{n_v n_w} \\
& & & 
\end{array}
\end{bmatrix}
\]
the components of the dimension vector $\alpha_e = (a_1,\hdots,a_{n_v},b_1,\hdots,b_{n_w})$ satisfy the following inequalities
\[
a_i \leq \sum_{j=1}^{n_w}  min(P_{ij},Q_{ij,v}) b_j \qquad \text{and} \qquad
b_j \leq \sum_{i=1}^{n_v}  min(P_{ij},R_{ij,v}) a_i \]
where $P_{ij}$ is the number of arrows in $\Gamma_e$ from vertex $\mu_i^{(v)}$ to vertex $\mu_j^{(w)}$ and where $Q_{ij,v}$ (resp. $R_{ij,v}$)  is the number of linearly independent entries in the $v$-th column (resp. in the $v$-th row) of the block-matrix $A_{ij}$.
\end{theorem}

In general, however,it is no longer true that the \'etale quiver associated to $\C \Gamma_{\Delta}$ will be symmetric.

\begin{example} (Continuing the $GL_2(\Z)$ example) Recall that the Zariski quiver $\Gamma(GL_2(\Z))$ is of the following form
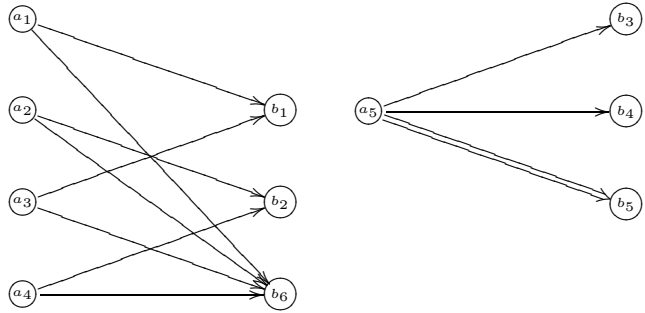
\[
\xymatrix{
\vtx{a_1} \ar[rrrd] \ar[rrrddd] & & & \\
\vtx{a_2} \ar[rrrd] \ar[rrrdd] & & & \vtx{b_1} \\
\vtx{a_3} \ar[rrru] \ar[rrrd] & & & \vtx{b_2} \\
\vtx{a_4} \ar[rrru] \ar[rrr] & & & \vtx{b_6} } \qquad
\xymatrix{
& & & \vtx{b_3} \\
\vtx{a_5} \ar[rrru] \ar[rrr] \ar@{=>}[rrrd] & & & \vtx{b_4} \\
& & & \vtx{b_5} } \]
Denote with $B_{ij}$ the uniquely determined arrow from the $i$-th left vertex to the $j$-th right vertex and with $B_{55}^{(1)}$ and $B_{55}^{(2)}$ the two arrows from $a_5$ to $b_5$. Consider the matrix
\[
\Delta = \begin{bmatrix}
B_{11} & 0 & B_{31} & 0 & 0 & 0  \\
0 & B_{22} & 0 & B_{42} & 0 & 0 \\
B_{16} & B_{26} & B_{36} & 0 & 0 & 0 \\
B_{16} & B_{26} & 0 & B_{46} & 0 & 0 \\
0 & 0 & 0 & 0 & B_{53} & 0 \\
0 & 0 & 0 & 0 & 0 & B_{54} \\
0 & 0 & 0 & 0 & B_{55}^{(1)} & B_{55}^{(2)} \\
0 & 0 & 0 & 0 & B_{55}^{(2)} & B_{55}^{(1)}
\end{bmatrix}
\]
Then, in addition to the generators $\{ g_1,\hdots, g_{10} \}$ described above, the semigroup $\mathcal{S}_{\Delta}$ has the additional generators
\[
\begin{array}{c|c|cccc|ccc}
& n & a_1 & a_2 & a_3 & a_4 & b_1 & b_2 & b_6 \\
\hline
g_{11} & 2 & 1 & 0 & 1 & 0 & 0 & 0 & 1 \\
g_{12} & 2 & 0 & 1 & 0 & 0 & 1 & 0 & 1 \\
\end{array}
\]
for the first component, and another two for the second component
\[
\begin{array}{c|c|c|ccc}
& n & a_5 & b_3 & b_4 & b_5 \\
\hline
g_{13} & 4 & 2 & 2 & 0 & 1 \\
g_{14} & 4 & 2 & 0 & 2 & 1 \\
\end{array}
\]
The \'etale quiver of $\C \Gamma_{\Delta}$ has the following form
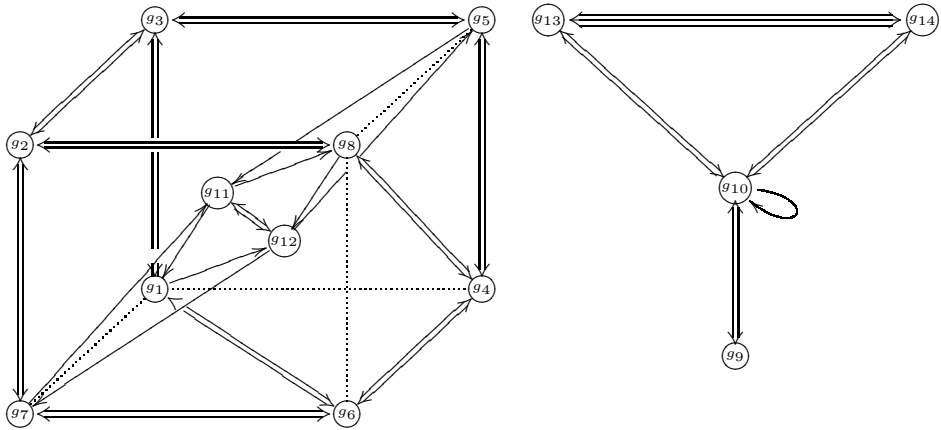
\[
\xymatrix@=.3cm{
& & & \vtx{g_3} \ar@{<=>}[rrrrrr] \ar@{<=>}'[ddd]'[dddddd]|(.77){\hole} \ar@{<=>}[lllddd] & & & & & & \vtx{g_5} \\
& & & & & & & & & \\
& & & & & & & & & \\
\vtx{g_2} \ar@{<=>}[rrrrrr] \ar@{<=>}[dddddd] & & & & & & \vtx{g_8} \ar@{.}[rrruuu] \ar@{.}[dddddd] \ar@{<=>}[rrrddd] & & & \\
& & & & \vtx{g_{11}}  \ar[ldd] \ar[rru] \ar@{<-}[rrrrruuuu]|(.27){\hole} \ar@{<-}[llllddddd] & & & & & \\
& & & & &  \vtx{g_{12}} \ar@{<=>}[lu] \ar@{<-}[lld] \ar@{<-}[ruu] \ar[rrrruuuuu]|(.35){\hole}& & & & \\
& & & \vtx{g_1} & & & & & & \vtx{g_4} \ar@{<=>}[uuuuuu] \ar@{.}[llllll] \\
& & & & & & & & & \\
& & & & & & & & & \\
\vtx{g_{7}} \ar@{<=>}[rrrrrr] \ar@{.}[rrruuu]  \ar@{<-}[rrrrruuuu] & & & & & & \vtx{g_6} \ar@{<=>}[llluuu]|(.88){\hole}  \ar@{<=>}[rrruuu] & & & } \quad
\xymatrix{
\vtx{g_{13}} & & & & \vtx{g_{14}} \ar@3@{<->}[llll] \\
& & & & \\
& & \vtx{g_{10}} \ar@{<=>}[lluu]  \ar@{<=>}[dd] \ar@{<=>}[rruu] \ar@(r,rd) & & \\
& & & & \\
& & \vtx{g_{9}} & & } \]
where the first component is not symmetric but has the interesting property that every vertex has exactly three incoming and three outgoing arrows.
\end{example}

\end{document}